\documentclass[12pt,twoside,reqno]{extarticle}
\usepackage[T1]{fontenc}
\usepackage[utf8]{inputenc}
\usepackage{float}
\usepackage{mathrsfs}
\usepackage{mathtools}
\usepackage{amsmath}
\usepackage{amsthm}
\usepackage{amssymb}
\usepackage{graphicx}
\usepackage{geometry}
\usepackage{microtype}
\usepackage[bookmarks=false,
 breaklinks=false,pdfborder={0 0 1},backref=false,colorlinks=false]
 {hyperref}
\hypersetup{pdftitle={Perimeter and exit-time variance in conformal Skorokhod embedding},
 pdfauthor={Maher Boudabra},
 pdfsubject={Perimeter comparisons, constant quantile allocation, and exit-time variance},
 pdfkeywords={Gross symmetrization, Brownian motion, perimeter, exit-time variance, Hardy spaces}}

\makeatletter
\@ifundefined{date}{}{\date{}}

\usepackage{mlmodern}

\usepackage[mathdesign=Garamond,TeX]{mathgreeks}
\mathgreeksfont{enc=OML,fam=mdugm,series=b,shape=n}
\input{omlmdugm.fd}
\DeclareFontShape{OML}{mdugm}{m}{it}
  {<-> s * [0.93] mdugmri7m}{}
\DeclareFontShape{OML}{mdugm}{m}{n}
  {<-> s * [0.93] mdugmr7m}{}

\usepackage{xcolor}
\usepackage{pdfrender}

\pdfrender{
  TextRenderingMode=2,
  LineWidth=0.1pt
}

\usepackage{tikz}

\usetikzlibrary{arrows.meta,positioning,fit,backgrounds,calc,shadows.blur}
\definecolor{ink}{HTML}{1F2A44}
\definecolor{degen}{HTML}{B23A48}
\definecolor{exist}{HTML}{2A6F97}
\definecolor{gross}{HTML}{2F7D4F}
\definecolor{spark}{HTML}{D98E04}
\tikzset{
  base/.style={rounded corners=3pt, draw=#1, line width=0.6pt, fill=#1!7,
     text width=5.1cm, align=center, inner sep=5pt, font=\small, text=ink,
     minimum height=1.25cm},
  head/.style={rounded corners=3pt, fill=#1, text=white, text width=5.6cm,
     align=center, inner sep=5pt, font=\small\bfseries, minimum height=0.9cm},
  flow/.style={-{Stealth[length=2.4mm]}, line width=0.7pt, color=ink!80},
  weak/.style={-{Stealth[length=2.2mm]}, dashed, line width=0.6pt, color=ink!55},
  tag/.style={font=\scriptsize\itshape, text=ink!70, fill=white, inner sep=1.5pt},
}

\usetikzlibrary{arrows.meta,calc}
\usetikzlibrary{arrows.meta}
\theoremstyle{plain}
\newtheorem{theorem}{Theorem}[section]
\newtheorem{proposition}[theorem]{Proposition}
\newtheorem{lemma}[theorem]{Lemma}
\newtheorem{corollary}[theorem]{Corollary}

\theoremstyle{definition}
\newtheorem*{example*}{Example}

\theoremstyle{remark}
\newtheorem{remark}[theorem]{Remark}

\numberwithin{equation}{section}

\newcommand{\D}{\mathbb D}
\newcommand{\T}{\mathbb T}
\newcommand{\C}{\mathbb C}
\newcommand{\R}{\mathbb R}

\newcommand{\Bsym}{\mathfrak B}
\newcommand{\cH}{\mathcal H}

\DeclareMathOperator{\Per}{Per}
\DeclareMathOperator{\Area}{Area}

\newcommand{\doi}[1]{%
  \href{https://doi.org/#1}{\nolinkurl{#1}}%
}

\makeatother

\begin{document}
\title{Beyond Area: Perimeter Optimization in the Planar Skorokhod Embedding Problem}

\author{Maher Boudabra%
\thanks{Department of Mathematics,
Monastir Preparatory Engineering Institute,
Monastir, Tunisia.
\newline
\href{mailto:maher.boudabra@gmail.com}
{\nolinkurl{maher.boudabra@gmail.com}}}}

\maketitle

\begin{abstract}
The planar Skorokhod embedding problem asks for a
simply connected planar domain whose Brownian exit
distribution has a prescribed real marginal $\mu$.
We investigate the minimization of conformal
perimeter over the class of all such embeddings.

We first establish the main existence theorem showing
that whenever a finite-perimeter embedding exists,
the infimum is attained (and vice versa). We derive necessary
conditions for finiteness and sufficient conditions
based on logarithmic integrability of the density.
In particular, we demonstrate that finite area
does not guarantee the existence of an embedding
with finite perimeter.

Although Gross's canonical embedding minimizes
area, we construct an explicit family of domains
for which Brownian symmetrization strictly
increases perimeter, disproving its universal
perimeter minimality. On the positive side, using
Pruss's containment theorem and inequalities
from convex geometry, we establish a general
lower bound in terms of the convex hull of the
Gross domain. Consequently, whenever the Gross
embedding is bounded and convex, it is the unique
perimeter minimizer.

Finally, we formulate an open problem concerning
the relationship between the existence of
finite-perimeter embeddings and the regularity
of Gross's canonical embedding.
\end{abstract}

\noindent\textbf{2020 Mathematics Subject Classification.}
Primary 49Q10, 60J65;
Secondary 30C70, 30H10, 60G40.

\medskip

\noindent\textbf{Keywords.}
Planar Skorokhod embedding,
conformal perimeter,
shape optimization,
Brownian motion,
Gross symmetrization,
Hardy spaces.

\section{Introduction and Problem Formulation}\label{Introduction}

Let $\left(Z_{t}\right)_{t\geq0}$ be a planar Brownian motion started
at $0$ with generator $\Delta/2$, and write 
\[
\tau_{U}=\inf\left\{ t\geq0:Z_{t}\notin U\right\} .
\]
In other words, $\tau_{U}$ is the exit time of $(Z)_{t \geq 0}$ from $U$. Now,
fix a real centered probability measure $\mu$ with a finite second
moment and look for a planar domain $U$ such that 
\[
\operatorname{Law}(\operatorname{Re}Z_{\tau_{U}})=\mu.
\]
This is the planar Skorokhod embedding problem PSEP, originally investigated
by Gross in \cite{Gross2019}. A solution is called $\mu$-domain
\cite{newMG}, or a planar embedding of $\mu$. In \cite{Gross2019},
Gross constructed a canonical solution: a simply connected domain,
which we denote always by $U^{G}_{\mu}$. Later on, Gross construction
was shown to be applicable to any $\mu$ having some finite $p$-th
moment with a $p>1$ \cite{mgremark} (the endpoint $p=1$ was discussed
in \cite{endpoint}). The geometric transformation of a domain $U\subseteq\mathbb{C}$
defined by 
\begin{equation}
\mathfrak{B}:U\longmapsto U^{G}_{\mu_{U}}\label{brown symm}
\end{equation}
where $\mu_{U}=\operatorname{Law}(\operatorname{Re}Z_{\tau_{U}})$
was introduced in \cite{Energy,Boudabra2026}. The map \ref{brown symm}
is referred to as Brownian symmetrization. This construction preserves
a marginal probability distribution determined by harmonic measure
at the origin. In \cite{Boudabra2026}, the author showed that the
Brownian approach provides a genuine extension of the classical Baernstein
symmetrization studied by Pruss in \cite{Pruss}, i.e. encompassing
a broader class of domains, namely unbounded ones. In the bounded
simply connected setting, however, the two constructions coincide.
The Baernstein symmetrization (attributed by Pruss to Baernstein)
rearranges the real boundary values of a conformal map symmetrically
and then reconstructs the imaginary part by harmonic conjugation,
normalized to vanish at the origin. If $f_{U}:\mathbb{D}\to U$ is
conformal and onto with $f(0)=0$, then $\mu_{U}$ is the distribution
of $\operatorname{Re}f^{*}(e^{i\Theta})$ for a uniform angle $\Theta$.
This identification connects the probabilistic construction with the
classical theory of symmetrization. Brownian symmetrization does not
increase the area. In other words 
\[
\operatorname{Area}(\mathfrak{B}(U))\leq\operatorname{Area}(U).
\]
In particular, Gross embedding has the least area. See \cite{Pruss,Boudabra2026}
for the two different perspectives.

Although the area minimization problem may also be
degenerate, its structure differs fundamentally
from that of conformal-perimeter minimization.
Indeed, Gross's embedding always minimizes area,
whereas it need not minimize conformal perimeter as it will be discussed.
Moreover, finite area does not guarantee the
existence of any finite-perimeter embedding:
there exist centered, compactly supported
probability measures for which the area minimization is nondegenerate while the perimeter minimization is degenerate.

Thus, perimeter optimization introduces both
a stronger regularity requirement on admissible
embeddings and a genuinely nontrivial
minimization problem. Our results address these
two aspects through existence and degeneracy
criteria, explicit counterexamples to Gross
minimality, and a convexity-based uniqueness
theorem.

Throughout, \textbf{$\mu$ is a centered, nondegenerate probability
measure on $\mathbb{R}$ with a bounded support}. In particular, $\mu$
has finite moments of any order. Denote by $\mathcal{C}(\mu)$ the
class of all proper simply connected domains $U\subsetneq\mathbb{C}$
such that $0\in U$ and 
\[
\operatorname{Law}(\operatorname{Re}Z_{\tau_{U}})=\mu,
\]
where, as usual, $(Z)_{t \geq 0}$ is a planar Brownian motion started at $0$.
An element of $\mathcal{C}(\mu)$ is called admissible embedding of
$\mu$. For each $U\in\mathcal{C}(\mu)$, write $f_{U}:\mathbb{D}\to U$
the onto normalized conformal map 
\[
f_{U}(0)=0,\qquad f_{U}'(0)>0.
\]
The conformal perimeter of $U$ is defined as 
\[
\operatorname{Per}_{\mathrm{c}}(U)=2\pi\|f_{U}'\|_{\mathbf{H}^{1}},
\]
where $\|\cdot\|_{\mathbf{H}^{1}}$ is the $1$-Hardy norm, with the
value $+\infty$ when $f_{U}'\notin\mathbf{H}^{1}$. In this manuscript,
we consider the minimization of the \emph{conformal perimeter} over
$\mathcal{C}(\mu)$ 
\begin{equation}
m_{\mathrm{c}}(\mu):=\inf_{U\in\mathcal{C}(\mu)}\operatorname{Per}_{\mathrm{c}}(U).\label{optimization problem}
\end{equation}
We shall call \ref{optimization problem} the conformal perimeter
optimization problem $\operatorname{CPOP}$. If $m_{\mathrm{c}}(\mu)$
is finite we say that $\operatorname{CPOP}$ is nondegenerate, otherwise
degenerate. The (increasing) quantile of $\mu$, which we denote by
$Q_{\mu}$, is defined by 
\[
Q_{\mu}(u)=\inf\left\{ x\in\mathbb{R}|\mu(-\infty,x]\geq u\right\} ,\,\,\,u\in(0,1).
\]
A famous statistical feature of $Q_{\mu}$ is that, when fed with
uniform values in $(0,1)$, it generates values sampling as $\mu$.
In other words, 
\[
Q_{\mu}(\operatorname{Uni}(0,1))\sim\mu.
\]
Set 
\[
q_{\mu}(t)=Q_{\mu}\left(1-\frac{|t|}{\pi}\right),\qquad-\pi<t<\pi,
\]
with a periodic extension and continuous endpoint values when available.
Its mean is zero. In this convention, the Gross map is 
\begin{equation}
g_{\mu}(z)=\sum_{n\ge1}c_{n}z^{n},\qquad c_{n}=\frac{2}{\pi}\int^{\pi}_{0}q_{\mu}(t)\cos(nt)\,dt.\label{eq:gross-series}
\end{equation}
The original choice of Gross is the increasing datum 
\[
Q_{\mu}\left(\frac{|t|}{\pi}\right),\qquad-\pi<t<\pi
\]
The two analytic completion maps differ by precomposition with $z\mapsto-z$
and have the same image. Gross's theorem gives the univalence of $g$,
$g(0)=0$, and $g(\D)=U^{G}_{\mu}$ is the Gross embedding.

To preserve the logical flow of the exposition, proofs containing
intermediate results that are used subsequently throughout the paper
are presented immediately after their statements. By contrast, purely
technical proofs serving only as auxiliary tools are deferred to the
appendix. When the symbol $\dagger$ is used next to a result, it
means the proof is deferred. Also, when $\mu$ is clear from the context,
it is omitted as a subscript.

\newpage{}

The results of this manuscript can be summarized as follows: \\

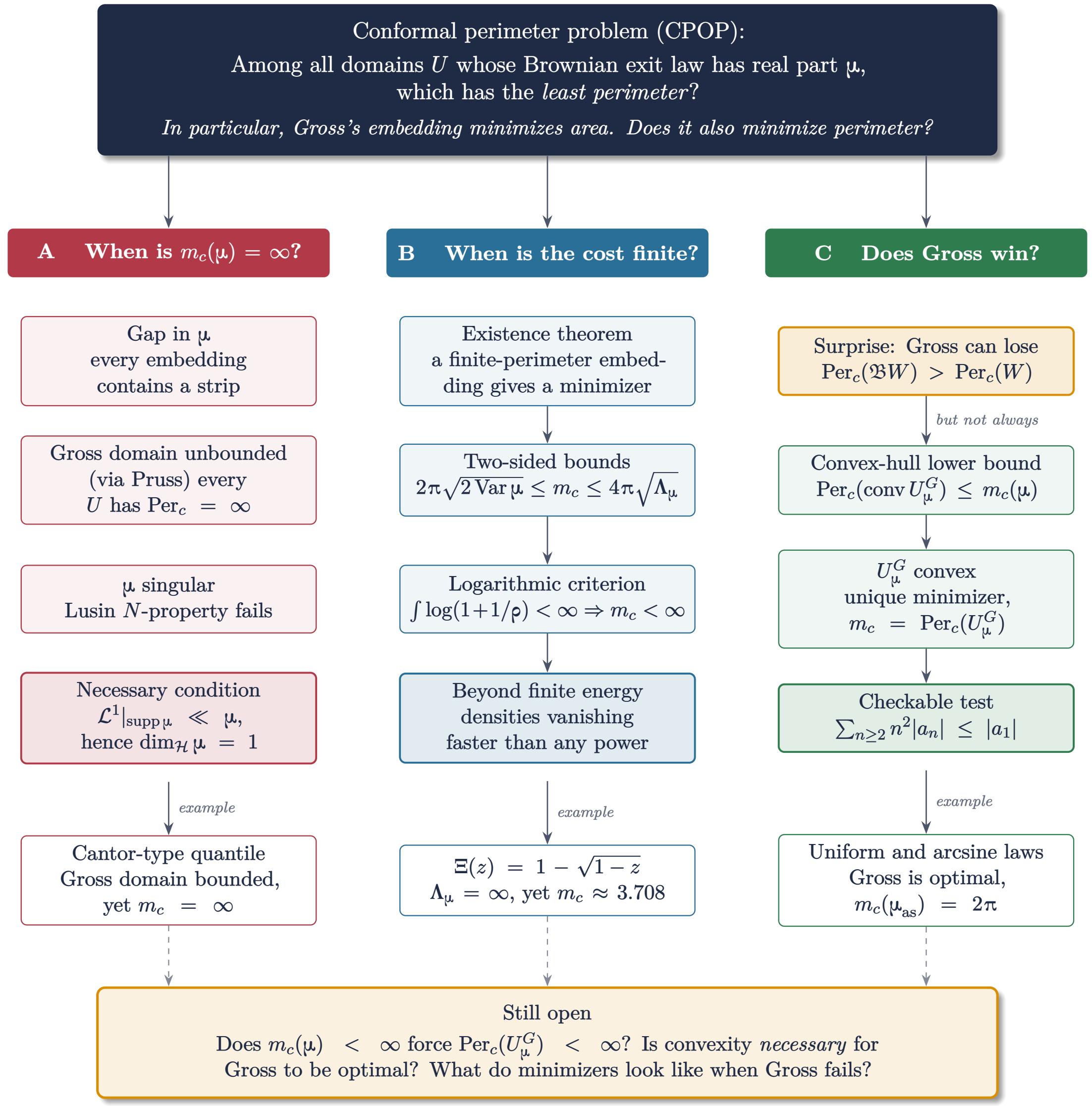
\begin{figure}[ht]
\centering
\resizebox{\textwidth}{!}{%
\begin{tikzpicture}[>=Stealth]

\node[rounded corners=4pt, fill=ink, text=white, text width=16cm, align=center,
      inner sep=9pt, font=\normalsize, blur shadow={shadow blur steps=5}] (Q) at (7,12.8)
  {Conformal perimeter problem (CPOP):\\[3pt]
   Among all domains $U$ whose Brownian exit law has real part $\mu$,\\
   which has the \emph{least perimeter}?\\[4pt]
   {\small\itshape In particular, Gross's embedding minimizes area. Does it also minimize perimeter?}};

\node[head=degen] (HA) at (0,9.6)  {A \;\; When is $m_c(\mu)=\infty$?};
\node[head=exist] (HB) at (7,9.6)  {B \;\; When is the cost finite?};
\node[head=gross] (HC) at (14,9.6) {C \;\; Does Gross win?};

\draw[flow] (Q.south -| HA) -- (HA.north);
\draw[flow] (Q.south -| HB) -- (HB.north);
\draw[flow] (Q.south -| HC) -- (HC.north);

\node[base=degen] (A1) at (0,7.6) {Gap in $\mu$\\ every embedding contains a strip};
\node[base=degen] (A2) at (0,5.4) {Gross domain unbounded\\ (via Pruss) every $U$ has $\mathrm{Per}_c=\infty$};
\node[base=degen] (A3) at (0,3.2) {$\mu$ singular\\ Lusin $N$-property fails};
\node[base=degen,fill=degen!14,line width=1pt] (A4) at (0,1.0)
  {Necessary condition\\ $\mathcal{L}^1|_{\operatorname{supp}\mu}\ll\mu$, hence $\dim_{\mathcal H}\mu=1$};

\node[base=exist] (B1) at (7,7.6) {Existence theorem\\ a finite-perimeter embedding gives a minimizer};
\node[base=exist] (B2) at (7,5.4) {Two-sided bounds\\ $2\pi\sqrt{2\operatorname{Var}\mu}\le m_c\le 4\pi\sqrt{\Lambda_\mu}$};
\node[base=exist] (B3) at (7,3.2) {Logarithmic criterion\\ $\int\log(1+1/\rho)<\infty\Rightarrow m_c<\infty$};
\node[base=exist,fill=exist!14,line width=1pt] (B4) at (7,1.0)
  {Beyond finite energy\\ densities vanishing faster than any power};
\draw[flow] (B1)--(B2); \draw[flow] (B2)--(B3); \draw[flow] (B3)--(B4);

\node[base=spark,fill=spark!16,line width=1.1pt] (C1) at (14,7.6)
  {Surprise: Gross can lose\\ $\mathrm{Per}_c(\mathfrak{B}W)>\mathrm{Per}_c(W)$};
\node[base=gross] (C2) at (14,5.4) {Convex-hull lower bound\\ $\mathrm{Per}_c(\operatorname{conv}U^G_\mu)\le m_c(\mu)$};
\node[base=gross] (C3) at (14,3.2) {$U^G_\mu$ convex\\ unique minimizer, $m_c=\mathrm{Per}_c(U^G_\mu)$};
\node[base=gross,fill=gross!14,line width=1pt] (C4) at (14,1.0) {Checkable test\\ $\sum_{n\ge2}n^2|a_n|\le|a_1|$};
\draw[flow] (C1)--(C2) node[midway,right=3pt,tag]{but not always};
\draw[flow] (C2)--(C3); \draw[flow] (C3)--(C4);

\begin{scope}[on background layer]
 \node[fit=(HA)(A1)(A4), draw=degen!70, rounded corners=6pt, fill=degen!3, inner sep=9pt] (FA) {};
 \node[fit=(HB)(B1)(B4), draw=exist!70, rounded corners=6pt, fill=exist!3, inner sep=9pt] (FB) {};
 \node[fit=(HC)(C1)(C4), draw=gross!70, rounded corners=6pt, fill=gross!3, inner sep=9pt] (FC) {};
\end{scope}

\node[base=degen,fill=white] (E1) at (0,-2.0)
  {Cantor-type quantile\\ Gross domain bounded, yet $m_c=\infty$};
\node[base=exist,fill=white] (E2) at (7,-2.0)
  {$\Xi(z)=1-\sqrt{1-z}$\\ $\Lambda_\mu=\infty$, yet $m_c\approx 3.708$};
\node[base=gross,fill=white] (E3) at (14,-2.0)
  {Uniform and arcsine laws\\ Gross is optimal, $m_c(\mu_{\mathrm{as}})=2\pi$};
\draw[flow] (FA.south -| E1) -- (E1.north) node[midway,right=3pt,tag]{example};
\draw[flow] (FB.south -| E2) -- (E2.north) node[midway,right=3pt,tag]{example};
\draw[flow] (FC.south -| E3) -- (E3.north) node[midway,right=3pt,tag]{example};

\node[rounded corners=4pt, draw=spark, line width=1.3pt, fill=spark!12, text=ink, text width=16cm,
      align=center, inner sep=9pt, font=\small, blur shadow={shadow blur steps=5}] (O) at (7,-5.0)
  {Still open\\[3pt]
   Does $m_c(\mu)<\infty$ force $\mathrm{Per}_c(U^G_\mu)<\infty$?
   Is convexity \emph{necessary} for Gross to be optimal?
   What do minimizers look like when Gross fails?};
\foreach \i in {1,2,3} \draw[weak] (E\i.south) -- (E\i.south |- O.north);

\end{tikzpicture}%
}

\caption{Overview of the main results and open
problems in conformal-perimeter optimization.}
\label{fig:research-flow}
\end{figure}

\section{Conformal and geometric perimeter }\label{sec:conformal-perimeter}

The notion of perimeter in the framework of simply connected domains,
although intuitive, needs a bit of carefulness, since we need to define
the notion of boundary length. We distinguish between the geometric
perimeter and the conformal perimeter. Let $1\le p<\infty$. For an
analytic function $f:\mathbb{D}\rightarrow\mathbb{C}$ , the $p$-Hardy
norm of $f$ is defined by 
\begin{equation}
\|f\|_{\mathbf{H}^{p}}=\sup_{0<r<1}\left(\frac{1}{2\pi}\int^{2\pi}_{0}|f(re^{it})|^{p}\,dt\right)^{1/p}.\label{eq:hardynorm}
\end{equation}
That is, the set of $f$ with finite $p$-Hardy norm is called $p$-Hardy
space, denoted by $\mathbf{H}^{p}$. A crucial property of Hardy spaces
is that when $f\in\mathbf{H}^{p}$ then $f$ admits a radial extension
$f^{*}$. More precisely, 
\[
f^{*}(\theta):=\lim_{r\uparrow1}f(re^{\theta i})
\]
exists a.e on the unit circle $\mathbb{T}$. Moreover, $f^{*}\in L^{p}(\mathbb{T})$.
We refer the reader to \cite{Duren1970,rudin} for detailed treatments
of the theory of Hardy spaces.

\subsection{Conformal perimeter}

For an analytic function $h$ on $\mathbb{D}$, write 
\[
M_{1}(r,h):=\frac{1}{2\pi}\int^{2\pi}_{0}|h(re^{it})|\,dt,\qquad0<r<1.
\]
The quantity $r\cdot M_{1}(r,h)$ represents the arc length of the
curve 
\[
\gamma_{r}(t)=h(re^{it}),\qquad0\leq t\leq2\pi.
\]
We write 
\[
L(\gamma_{r})=2\pi rM_{1}(r,h').
\]

Let $f$ a normalized univalent function mapping $\mathbb{D}$ onto
$U$. We define the \emph{conformal perimeter} of $U$ by 
\begin{equation}
\Per_{\mathrm{c}}(U):=2\pi\|f'\|_{\mathbf{H}^{1}},\label{eq:conformal-perimeter}
\end{equation}
with the convention that $\Per_{\mathrm{c}}(U)=+\infty$ when $f'\notin\mathbf{H}^{1}$.
This definition is independent of the choice of normalized conformal
map. Indeed, any two conformal maps from $\mathbb{D}$ onto $U$ differ
by precomposition with an automorphism of the unit disk, and the boundary
length represented by \eqref{eq:conformal-perimeter} is invariant
under such a change of parametrization. Since the membership of $f'$
in $\mathbf{H}^{1}$ implies that of $f$, then by the radial extension
in Hardy spaces, one has this direct geometric interpretation: $\Per_{\mathrm{c}}(U)$
is the length of the parametrized curve
\begin{equation}
\xi\longmapsto(f^{*}(\xi)\label{eq:trace-length}
\end{equation}
where $f^{*}$ is the radial limit of $f$. Thus the conformal perimeter
is the length of the boundary as traversed by the conformal boundary
parametrization, with any multiplicity of traversal retained. Note
that if $\Per_{\mathrm{c}}(U)$ is finite then $f$ extends continuously
to $\overline{\mathbb{D}}$. However $f$ may loose the injectivity
in the closed disc.

\subsection{Geometric perimeter}

For a set $E\subset\mathbb{C}$, its one-dimensional Hausdorff measure
is defined by

\[
\mathcal{H}^{1}(E)=\lim_{\delta\downarrow0}\mathcal{H}^{1}_{\delta}(E),
\]
where

\[
\mathcal{H}^{1}_{\delta}(E):=\inf\left\{ \sum^{\infty}_{j=1}\operatorname{diam}(E_{j}):E\subset\bigcup^{\infty}_{j=1}E_{j},\quad\operatorname{diam}(E_{j})\le\delta\right\} .
\]
Accordingly, for a simply connected domain $U$ we define the \emph{geometric
perimeter} of $U$ by

\begin{equation}
\Per_{\mathrm{g}}(U):=\mathcal{H}^{1}(\partial U)\label{eq:geometric-perimeter}
\end{equation}
whenever this quantity is finite. Geometrically, $\mathcal{H}^{1}(\partial U)$
is the total one-dimensional length of the boundary viewed simply
as a subset of the plane. It counts each geometric boundary point
only once. So it is the natural generalization of ordinary curve length
to boundaries that may have corners, cusps, slits, or more complicated
shapes. A crucial distinction is that 
\[
\Per_{\mathrm{g}}(U)\leq\Per_{\mathrm{c}}(U)
\]
in general. Consider the slit domain 
\[
U=\mathbb{D}\setminus[1/2,1).
\]
It is simply connected. We can show that 
\[
\Per_{\mathrm{g}}(U)=2\pi+\frac{1}{2}<2\pi\|f'\|_{\mathbf{H}^{1}}=2\pi+2\left(\frac{1}{2}\right).
\]
The strict inequality comes from the fact that Hausdorff measure counts
the slit only once, whereas the conformal boundary parametrization
sees its two sides separately.

\subsection{Relation between the two notions}

For a rectifiable simple closed curve, $\mathcal{H}^{1}$ agrees with
its ordinary arc length. Hence, if $U$ is a rectifiable Jordan domain,
$\Per_{\mathrm{g}}(U)$ is the usual Euclidean perimeter. The distinction
between \eqref{eq:conformal-perimeter} and \eqref{eq:geometric-perimeter}
is relevant when the boundary is not Jordan. The conformal boundary
parametrization may traverse the same geometric boundary piece more
than once. In that case the conformal perimeter records the multiplicity,
whereas Hausdorff measure counts each geometric point only once. The
following implications summarize the situation.

\begin{proposition}\label{prop:conformal-geometric-perimeter} \cite{Bishop}
Let $U\subsetneq\mathbb{C}$ be simply connected and let $f:\mathbb{D}\to U$
be conformal with range $U$. We have 
\begin{equation}
\Per_{\mathrm{g}}(U)\le\Per_{\mathrm{c}}(U)\le2\Per_{\mathrm{g}}(U).\label{per geo comparable}
\end{equation}
In particular, both geometric and conformal perimeters share the same
finiteness status. If, in addition, $U$ is a Jordan domain, then
\[
f'\in\mathbf{H}^{1}\quad\Longleftrightarrow\quad\partial U\text{ is rectifiable},
\]
and in this case 
\begin{equation}
\Per_{\mathrm{g}}(U)=\Per_{\mathrm{c}}(U)=2\pi\|f'\|_{\mathbf{H}^{1}}.\label{eq:perimeters-equal}
\end{equation}

\end{proposition}

The estimate \ref{per geo comparable} is referred to as Bishop's
inequality. The equivalence between the rectifiability of $\partial U$
and the membership of $f'$ in $\mathbf{H}^{1}$ is standard. The
factor $2$ in Bishop's inequality reflects the fact that the conformal
boundary parametrization traverses almost every boundary point at
most twice. Higher multiplicities can occur only on an exceptional
countable set and therefore do not contribute to boundary length.
Note that a simply connected domain obtained by removing a slit may
have a conformal map with derivative in $\mathbf{H}^{1}$, while the
boundary is not Jordan. The two sides of the slit correspond to distinct
boundary approaches and are traversed separately by the conformal
boundary parametrization. Consequently the slit is counted twice by
$\Per_{\mathrm{c}}$ but only once by $\Per_{\mathrm{g}}$.

\subsection{Pruss's containment theorem}

For a bounded open set $E\subset\C$, define the length of its vertical
section at $x\in\R$ by 
\[
\ell_{E}(x)=\mathcal{L}^{1}\bigl(\{y\in\R:x+iy\in E\}\bigr),
\]
where $\mathcal{L}^{1}$ denotes one-dimensional Lebesgue measure.
Its Steiner symmetrization about the real axis is 
\[
\mathscr{S}(E)=\left\{ x+iy:\ell_{E}(x)>0,\quad|y|<\frac{\ell_{E}(x)}{2}\right\} .
\]

Thus Steiner symmetrization centers each vertical section while preserving
its length. Steiner symmetrization is known to preserve the area while
not boosting the perimeter \cite{kesevan}. In \cite{Pruss}, Pruss
established the following containment theorem, which will play a role
to identify when Gross embedding is the minimizer of our perimeter
optimization problem.

\begin{theorem} \label{Pruss's cont} Let $U$ be a bounded simply
connected domain. Then 
\[
\mathfrak{B}(U)\subset\mathscr{S}(U).
\]

\end{theorem}

Note that the original statement of Pruss's containment theorem \ref{Pruss's cont}
is within a Hardy-space framework, which is automatic here because
$U$ is bounded. For a compact convex body $K$, we use the corresponding
closed symmetral 
\[
\mathscr{S}(K)=\left\{ x+iy:x\in\operatorname{proj}_{\R}K,\quad|y|\leq\frac{\ell_{K}(x)}{2}\right\} .
\]
This is again a compact convex body. Indeed, the upper endpoint of
a vertical section of $K$ is a concave function of $x$, and the
lower endpoint is convex, so their difference $\ell_{K}$ is concave.

We emphasize that \ref{Pruss's cont} does not, by itself, imply a
perimeter inequality: inclusion does not order the perimeters of arbitrary
nonconvex domains.

\section{Existence of a conformal-perimeter minimizer}\label{sec:perimeter-existence}

In this section, we establish our main existence theorem for the conformal-perimeter
minimization problem. We then investigate necessary and sufficient
conditions for the infimum to be finite, identifying both geometric
and measure-theoretic criteria.

\subsection{Existence theorem}

\begin{theorem}[Attainment of the conformal-perimeter infimum] \label{thm:perimeter-existence}
Suppose that there exists $U_{0}\in\mathcal{C}(\mu)$ such that 
\[
\operatorname{Per}_{\mathrm{c}}(U_{0})<\infty.
\]
Then the conformal-perimeter infimum is finite and attained. More
precisely, there exists $U_{*}\in\mathcal{C}(\mu)$ satisfying 
\begin{equation}
\operatorname{Per}_{\mathrm{c}}(U_{*})=\min_{U\in\mathcal{C}(\mu)}\operatorname{Per}_{\mathrm{c}}(U)<\infty.\label{eq:perimeter-attainment}
\end{equation}
The minimizing domain $U_{*}$ is bounded and its normalized conformal
map satisfies $f_{U_{*}}'\in\mathbf{H}^{1}$. No Jordan-boundary assumption
is imposed or concluded.

\end{theorem} 
\begin{proof}
We divide the argument into four steps.

\medskip{}
 \noindent \textit{Step 1. Uniform coefficient bounds.} By hypothesis,
\[
m_{\mathrm{c}}(\mu)\le\operatorname{Per}_{\mathrm{c}}(U_{0})<\infty.
\]
Choose a minimizing sequence $(U_{j})$ and let 
\[
f_{j}(z)=\sum^{\infty}_{n=1}c_{j,n}z^{n}
\]
be the corresponding normalized conformal maps. After discarding finitely
many terms, there exists $C<\infty$ such that 
\begin{equation}
\|f_{j}'\|_{\mathbf{H}^{1}}\le C\qquad(j\ge1).\label{eq:uniform-derivative-H1}
\end{equation}
For every $0<r<1$, Cauchy's coefficient formula applied to $f_{j}'$
gives 
\[
nc_{j,n}r^{n-1}=\frac{1}{2\pi}\int^{2\pi}_{0}f_{j}'(re^{it})e^{-i(n-1)t}\,dt.
\]
Therefore 
\[
n|c_{j,n}|r^{n-1}\le\frac{1}{2\pi}\int^{2\pi}_{0}|f_{j}'(re^{it})|\,dt\le C.
\]
Letting $r\uparrow1$, we obtain 
\begin{equation}
|c_{j,n}|\le\frac{C}{n},\qquad j,n\ge1.\label{eq:perimeter-coefficient-bound}
\end{equation}
In particular, 
\begin{equation}
\sum_{n>N}|c_{j,n}|^{2}\le C^{2}\sum_{n>N}\frac{1}{n^{2}}\le\frac{C^{2}}{N}.\label{eq:perimeter-H2-tails}
\end{equation}
Thus the Taylor coefficients have uniformly small $\mathbf{H}^{2}$
tails.

\medskip{}
 \noindent \textit{Step 2. Strong $\mathbf{H}^{2}$ compactness.}
For every fixed $n$, the sequence $(c_{j,n})_{j}$ is bounded. By
a diagonal argument, we may extract a subsequence such that 
\[
c_{j,n}\longrightarrow c_{n}\qquad\text{for every }n\ge1.
\]
Set 
\[
f(z)=\sum^{\infty}_{n=1}c_{n}z^{n}.
\]
The estimate $|c_{n}|\le C/n$ shows that this series defines an analytic
function in $\mathbb{D}$. Moreover, convergence of individual coefficients,
together with \eqref{eq:perimeter-coefficient-bound}, gives 
\[
f_{j}\longrightarrow f\qquad\text{locally uniformly in }\mathbb{D}.
\]
We claim that the convergence is in fact strong in $\mathbf{H}^{2}$.
Indeed, 
\[
|c_{j,n}-c_{n}|^{2}\le\frac{4C^{2}}{n^{2}},
\]
and the right-hand side is summable. Dominated convergence for series
therefore gives 
\begin{equation}
\|f_{j}-f\|^{2}_{\mathbf{H}^{2}}=\sum^{\infty}_{n=1}|c_{j,n}-c_{n}|^{2}\longrightarrow0.\label{eq:perimeter-strong-H2}
\end{equation}

\medskip{}
 \noindent \textit{Step 3. Univalence and preservation of the exit
law.} By the classical Hardy-space boundary regularity theorem, the
finite-perimeter competitor $U_{0}$ is bounded. Since it realizes
$\mu$, the latter has bounded support. In particular, 
\[
v:=\int_{\mathbb{R}}x^{2}\,\mu(dx)
\]
is finite and strictly positive. Let $m$ denote normalized angular
measure on $\mathbb{T}$, and put 
\[
u_{j}=\operatorname{Re}f^{*}_{j}.
\]
Conformal invariance of harmonic measure gives 
\[
(u_{j})_{*}m=\mu\qquad\text{for every }j.
\]
Since $f_{j}(0)=0$, Parseval's identity yields 
\[
\|f_{j}\|^{2}_{\mathbf{H}^{2}}=2\|u_{j}\|^{2}_{L^{2}(\mathbb{T},m)}=2v.
\]
By \eqref{eq:perimeter-strong-H2}, 
\[
\|f\|^{2}_{\mathbf{H}^{2}}=2v>0.
\]
Thus $f$ is nonconstant. The locally uniform limit of univalent analytic
functions is either univalent or constant. Consequently, $f$ is univalent.
Furthermore, 
\[
f(0)=0,\qquad f'(0)=\lim_{j\to\infty}f_{j}'(0)>0,
\]
where strict positivity follows from univalence of the limit. Next,
strong $\mathbf{H}^{2}$ convergence implies 
\[
\operatorname{Re}f^{*}_{j}\longrightarrow\operatorname{Re}f^{*}\qquad\text{in }L^{2}(\mathbb{T},m).
\]
In particular, these boundary functions converge in probability. Since
all of them have the same law $\mu$, their limit has that law as
well: 
\begin{equation}
(\operatorname{Re}f^{*})_{*}m=\mu.\label{eq:perimeter-law-preservation}
\end{equation}

\medskip{}
 \noindent \textit{Step 4. Lower semicontinuity and attainment.}
Local uniform convergence of the analytic maps implies 
\[
f_{j}'\longrightarrow f'\qquad\text{locally uniformly in }\mathbb{D}.
\]
For every fixed $0<r<1$, 
\[
\frac{1}{2\pi}\int^{2\pi}_{0}|f'(re^{it})|\,dt=\lim_{j\to\infty}\frac{1}{2\pi}\int^{2\pi}_{0}|f_{j}'(re^{it})|\,dt.
\]
Hence 
\[
\frac{1}{2\pi}\int^{2\pi}_{0}|f'(re^{it})|\,dt\le\liminf_{j\to\infty}\|f_{j}'\|_{\mathbf{H}^{1}}.
\]
Taking the supremum over $r$ gives 
\begin{equation}
\|f'\|_{\mathbf{H}^{1}}\le\liminf_{j\to\infty}\|f_{j}'\|_{\mathbf{H}^{1}}.\label{eq:perimeter-lsc}
\end{equation}
In particular, $f'\in\mathbf{H}^{1}$. Again , $f$ extends continuously
to $\overline{\mathbb{D}}$ and its image 
\[
U_{*}:=f(\mathbb{D})
\]
is bounded. Conformal invariance of Brownian motion, together with
\eqref{eq:perimeter-law-preservation}, shows that 
\[
\operatorname{Law}\bigl(\operatorname{Re}Z_{\tau_{U_{*}}}\bigr)=\mu.
\]
Thus $U_{*}\in\mathcal{C}(\mu)$. Finally, \eqref{eq:perimeter-lsc}
gives 
\[
\begin{aligned}\operatorname{Per}_{\mathrm{c}}(U_{*}) & =2\pi\|f'\|_{\mathbf{H}^{1}}\\
 & \le\liminf_{j\to\infty}2\pi\|f_{j}'\|_{\mathbf{H}^{1}}\\
 & =m_{\mathrm{c}}(\mu).
\end{aligned}
\]
The reverse inequality follows from the definition of $m_{\mathrm{c}}(\mu)$,
since $U_{*}$ is admissible. Therefore 
\[
\operatorname{Per}_{\mathrm{c}}(U_{*})=m_{\mathrm{c}}(\mu),
\]
which proves \eqref{eq:perimeter-attainment}. 
\end{proof}

\begin{corollary} One has 
\[
\exists\,U\in\mathcal{C}(\mu):\operatorname{Per}_{\mathrm{c}}(U)<\infty\,\,\Longleftrightarrow\,\,m_{\mathrm{c}}(\mu)<\infty
\]
and 
\[
\min_{U\in\mathcal{C}(\mu)}\operatorname{Per}_{\mathrm{c}}(U)=\inf_{U\in\mathcal{C}(\mu)}\operatorname{Per}_{\mathrm{c}}(U).
\]

\end{corollary}

Theorem~\ref{thm:perimeter-existence} does not require convexity,
symmetry, or a prescribed number of boundary branches. In particular,
it applies to the full class of planar embeddings of $\mu$. The existence
theorem \ref{thm:perimeter-existence} requires the existence of at
least a competitor. So far, only two embeddings are available: Gross
and Maher-Greg embeddings. As Maher-Greg embedding is always unbounded
\cite{newMG,Boudabra2026}, then it can never have a finite conformal
perimeter (it will be shown) . Among these two canonical constructions,
only the Gross embedding can provide a finite-perimeter competitor
for the $\operatorname{CPOP}$. This does not exclude finite-perimeter
embeddings obtained by other constructions. 

\begin{corollary}\label{cor:gross-finite-perimeter-existence} If
\[
\operatorname{Per}_{\mathrm{c}}(U^{G}_{\mu})<\infty
\]
then there exists an embedding $U_{*}$ of $\mu$ such that 
\[
\operatorname{Per}_{\mathrm{c}}(U_{*})=\min_{U\in\mathcal{C}(\mu)}\operatorname{Per}_{\mathrm{c}}(U).
\]
\end{corollary}

\begin{remark}[Scope of the existence result] The minimizing domain
need not be the Gross domain as it will be shown in the sequel. Moreover,
the theorem concerns conformal perimeter, which counts boundary traversal
with multiplicity. It does not establish that the minimizing domain
has a Jordan boundary, nor does it automatically give attainment of
the geometric-perimeter infimum within the restricted class of rectifiable
Jordan domains. \end{remark}

\subsection{Necessary conditions}

In this section, we gather some conditions that must hold if the $\operatorname{CPOP}$
is nondegenerate. Hence, these conditions serve as a primary test
to check before starting looking for the minimizer of the $\operatorname{CPOP}$.
Recall that in section \ref{Introduction}, we assumed that $\mu$
has bounded support. It is not a technical assumption, rather a necessity
to avoid the degeneracy of the $\operatorname{CPOP}$. A distribution
with unbounded support generates always embeddings with infinite perimeters.
This is based on the following lemma.

\begin{lemma} \label{bishopestimate} \cite{Bishop2} Let $U$ be
a simply connected domain. Then 
\begin{equation}
2\operatorname{diam}(U)\leq\operatorname{Per}_{\mathrm{c}}(U).\label{bishop estimate}
\end{equation}
In particular, if $U$ has a finite conformal perimeter then it is
bounded. \end{lemma}

Lemma \ref{bishopestimate} implies that if our $\operatorname{CPOP}$
is nondegenrate then $\mu$ has a bounded support. Moreover, from
\ref{bishop estimate}, we obtain 
\begin{equation}
m_{c}(\mu)\ge2\operatorname{diam}(\operatorname{supp}\mu).\label{supportbound}
\end{equation}
The support bound in \ref{supportbound} is sharp as it is asymptotically
attained by increasingly thin rectangular domains. We reiterate that
$\mu$ is always with a bounded support.

\begin{remark} Note that the bounded support of $\mu$ does not prevent
from obtaining unbounded embeddings. \end{remark}

The classical isoperimetric inequality asserts that 
\[
4\pi\operatorname{Area}(U)\leq\operatorname{Per}_{g}(U)^{2}
\]
for every planar set of finite area and finite geometric perimeter.
In our conformal setting, no boundedness assumption is needed. Indeed,
finite conformal perimeter implies that the normalized conformal map
has derivative in $\mathbf{H}^{1}$, and hence extends continuously
to the closed unit disk. In particular, the domain is bounded and
has finite area. Consequently, 
\begin{equation}
4\max(\operatorname{diam}(U)^{2},\pi\operatorname{Area}(U))\leq\operatorname{Per}_{\mathrm{c}}(U)^{2}.\label{isoperineq}
\end{equation}
Note that $\operatorname{diam}(U)^{2}$ and $\pi\operatorname{Area}(U)$
are uncorrelated.

\begin{proposition} If $\mu$ has a genuine gap $(a,b)$, i.e. $\mu((a,b))=0$
for some $a<b$, with positive mass on both sides of this interval,
then any planar embedding of $\mu$ contains the vertical strip $(a,b)+i\mathbb{R}$.
\end{proposition} 
\begin{proof}
Let $S=(a,b)+i\mathbb{R}$. Suppose, for a contradiction, that there
exists $\xi\in\partial U$ such that 
\[
a<\operatorname{Re}\xi<b.
\]
Choose $\varepsilon>0$ sufficiently small that $\operatorname{Ball}(\xi,\varepsilon)\subset S$.
By the classical full-support property of harmonic measure, for which
Markowsky gave a simple probabilistic proof \cite{markowskydirichlet},
followed by the strong maximum principle, we have 
\[
\mathbf{P}_{0}\bigl(Z_{\tau_{U}}\in\operatorname{Ball}(\xi,\varepsilon)\bigr)>0.
\]
Consequently, 
\[
\begin{aligned}\mu((a,b)) & =\mathbf{P}_{0}\bigl(a<\operatorname{Re}Z_{\tau_{U}}<b\bigr)\\
 & \ge\mathbf{P}_{0}\bigl(Z_{\tau_{U}}\in\operatorname{Ball}(\xi,\varepsilon)\bigr)>0,
\end{aligned}
\]
contradicting the assumption $\mu((a,b))=0$. Hence 
\[
\partial U\cap S=\varnothing.
\]
Since $S$ is connected and does not meet $\partial U$, we have either
$S\subset U$ or $S\cap U=\varnothing$. The latter is impossible:
by connectedness, it would force $U$ to lie entirely on one side
of the strip, contradicting the assumption that $\mu$ assigns positive
mass to both $(-\infty,a]$ and $[b,\infty)$. Therefore, 
\[
(a,b)+i\mathbb{R}\subset U.
\]
\end{proof}

\begin{corollary} If $\mu$ has a genuine gap $(a,b)$ then the $\operatorname{CPOP}$
is degenerate. \end{corollary}

\begin{theorem}[Unbounded Gross domains] \label{thm:unbounded-gross}
If the Gross embedding $U^{G}_{\mu}$ is unbounded then 
\begin{equation}
\operatorname{Per}_{\mathrm{c}}(U)=+\infty\qquad\text{for every }U\in\mathcal{C}(\mu).\label{eq:infinite-conformal-perimeter}
\end{equation}
In particular, the $\operatorname{CPOP}$ is degenerate. Equivalently,
\[
\exists\,U\in\mathcal{C}(\mu):\operatorname{Per}_{\mathrm{c}}(U)<\infty\quad\Longrightarrow\quad U^{G}_{\mu}\text{ is bounded}.
\]

\end{theorem} 
\begin{proof}
Suppose that there exists an embedding $U\in\mathscr{C}(\mu)$ such
that $\operatorname{Per}_{\mathrm{c}}(U)<\infty$. In particular,
$U$ is bounded by the isperimetric inequality \ref{isoperineq}.
By Pruss's containment theorem \ref{Pruss's cont}, 
\[
\mathfrak{B}(U)=U^{G}_{\mu}\subset\mathscr{S}(U),
\]
Since $U$ is bounded, its Steiner symmetrization $\mathscr{S}(U)$
is also bounded. Hence $U^{G}_{\mu}$ is bounded as well, contrary
to our assumption. Therefore every embedding of $\mu$ is unbounded.
Again by the isoperimetric inequality, we get 
\[
\operatorname{Per}_{\mathrm{c}}(U)=+\infty\qquad\text{for every }U\in\mathcal{C}(\mu).
\]
\end{proof}

\textbf{Comment.} The preceding result carried by theorem \ref{thm:unbounded-gross}
reveals a distinctive role of the Gross embedding in the perimeter
optimization problem. Its unboundedness is not merely a geometric
feature of a particular construction, but an intrinsic obstruction
imposed by the prescribed exit law. This property is specific to the
canonical role of the Gross embedding: the unboundedness of an arbitrary
admissible domain carries no such implication. For example, Maher-Greg
embedding is always unbounded, however it does not mean that the underlying
$\operatorname{CPOP}$ is degenerate.\\
 \\
 \\
 One might think that Gross embedding would be bounded if we pick
a distribution $\mu$ with no gaps, along with absolute continuity.
Surprisingly, this does not prevent the unboundedness of Gross embedding.
In their study of the boundedness of Gross domains, Boudabra, Kroumi,
and Smii \cite{BKS2024} established a characterization of the boundedness
of Gross embedding when the underlying distribution is continuous
and compactly supported. In particular, the following distribution
$\mu$ with density 
\[
\rho(x)=\left(2\int^{1}_{0}e^{-1/t}\,dt\right)^{-1}\exp\left(-\dfrac{1}{|x|}\right)\mathbf{1}_{(-1,1)}(x)
\]
generates an unbounded Gross domain. Hence, the associated $\operatorname{CPOP}$
is degenerate. Thus even for an absolutely continuous, compactly supported
law with no gaps, the conformal-perimeter optimization problem may
have no finite-valued solution. This example demonstrates that the
finite-perimeter competitor hypothesis in Theorem \ref{thm:perimeter-existence}
is a genuine restriction and cannot be replaced merely by regularity
of the distribution or connectedness of its support.\\
 \\

We now turn to measure-theoretic obstructions to the existence of
finite-perimeter embeddings. While boundedness of the Gross domain
provides a necessary geometric condition, the regularity of the prescribed
exit law imposes further restrictions. We first show that purely singular
distributions cannot admit finite-perimeter embeddings. We then establish
a stronger necessary condition: the absolutely continuous component
of the exit law must have a strictly positive density almost everywhere
on its support, although a singular component may still be present.

Two measures $\mu$ and $\nu$ are mutually singular, denoted by $\mu\perp\nu$,
if there exists a measurable set $E$ such that 
\[
\mu(E)=0\qquad\text{and}\qquad\nu(E^{c})=0.
\]
In other words, the two measures are concentrated on disjoint measurable
sets. The notation $\nu\ll\mu$ means that $\nu$ is absolutely continuous
with respect to $\mu$, namely, 
\[
\mu(E)=0\quad\Longrightarrow\quad\nu(E)=0
\]
for every measurable set $E$. In other words, every $\mu$-null set
is also $\nu$-null. Recall that, by the Lebesgue decomposition theorem
(see \cite{folland}), every Borel probability measure $\mu$ on $\mathbb{R}$
admits a unique decomposition 
\[
\mu(dx)=\rho(x)\,dx+\mu_{s}(dx),
\]
where $\rho\geq0$ is integrable and $\mu_{s}$ is singular with respect
to Lebesgue measure.

\begin{theorem} \label{thm:singular-law} Let $\mu$ be singular
with respect to Lebesgue measure. Then $U\in\mathcal{C}(\mu)$ satisfies
\[
\operatorname{Per}_{\mathrm{c}}(U)=+\infty.
\]
Consequently, the $\operatorname{CPOP}$ is degenerate\end{theorem} 
\begin{proof}
Suppose, by contradiction, that there exists $U\in\mathcal{C}(\mu)$
with finite conformal perimeter. In particular, we have $f_{U}'\in\mathbf{H}^{1}$
and hence $f_{U}$ extends continuously to $\overline{\mathbb{D}}$,
and its boundary trace is absolutely continuous. In particular, 
\[
u(\theta)=\operatorname{Re}f_{U}(e^{i\theta})
\]
is a nonconstant absolutely continuous function on $[0,2\pi]$. Since
$\mu$ is singular, there exists a Borel set $N\subset\operatorname{supp}(\mu)$
such that 
\[
\mathcal{L}^{1}(N)=0,\qquad\mu(N)=1.
\]
Because $\mu$ is the distribution of $u(\Theta)$ , where $\Theta$
is uniformly distributed on $[0,2\pi]$, the set 
\[
E=\{\theta:u(\theta)\notin N\}
\]
has Lebesgue measure zero. Absolute continuity implies Lusin's $N$-property:
the image under $u$ of every Lebesgue-null set has Lebesgue measure
zero (see \cite{Leoni}). so $\mathcal{L}^{1}(u(E))=0$. Moreover,
\[
u([0,2\pi])\subset N\cup u(E).
\]
It follows that $u([0,2\pi])$ has Lebesgue measure zero. This is
impossible, since $u$ is continuous and nonconstant, and therefore
its image is a nondegenerate interval. Thus no embedding of $\mu$
can have finite geometric perimeter. The equivalence of perimeter
finiteness gives the corresponding conclusion for conformal perimeter. 
\end{proof}

\begin{example*}
Let $(\xi_{n})_{n\ge1}$ be independent Bernoulli random variables
with 
\[
\mathbf{P}(\xi_{n}=1)=\frac{1}{4},\qquad\mathbf{P}(\xi_{n}=0)=\frac{3}{4},
\]
and define 
\[
X=\sum^{\infty}_{n=1}\frac{\xi_{n}}{2^{n}},\qquad\mu=\operatorname{Law}\left(X-\frac{1}{4}\right).
\]
$X$ is well defined and $\mathbf{E}(X)=1/4$. The measure $\mu$
is then centered. Moreover, $\mu$ is nonatomic, singular continuous,
and has support $[-1/4,3/4]$. Indeed, every dyadic interval has positive
probability, so the support of $X$ is $[0,1]$. On the other hand,
by the strong law of large numbers, almost every realization of $X$
has binary digit frequency $1/4$, whereas Lebesgue-almost every point
of $[0,1]$ has binary digit frequency $1/2$. Hence $\mu$ is singular
with respect to Lebesgue measure. Hence the underlying $\operatorname{CPOP}$
is degenerate by theorem \ref{thm:singular-law}. Note that we can
show that the Gross embedding is nevertheless bounded. 
\end{example*}
\begin{figure}[H]
\begin{centering}
\includegraphics[width=9cm,totalheight=9cm,keepaspectratio]{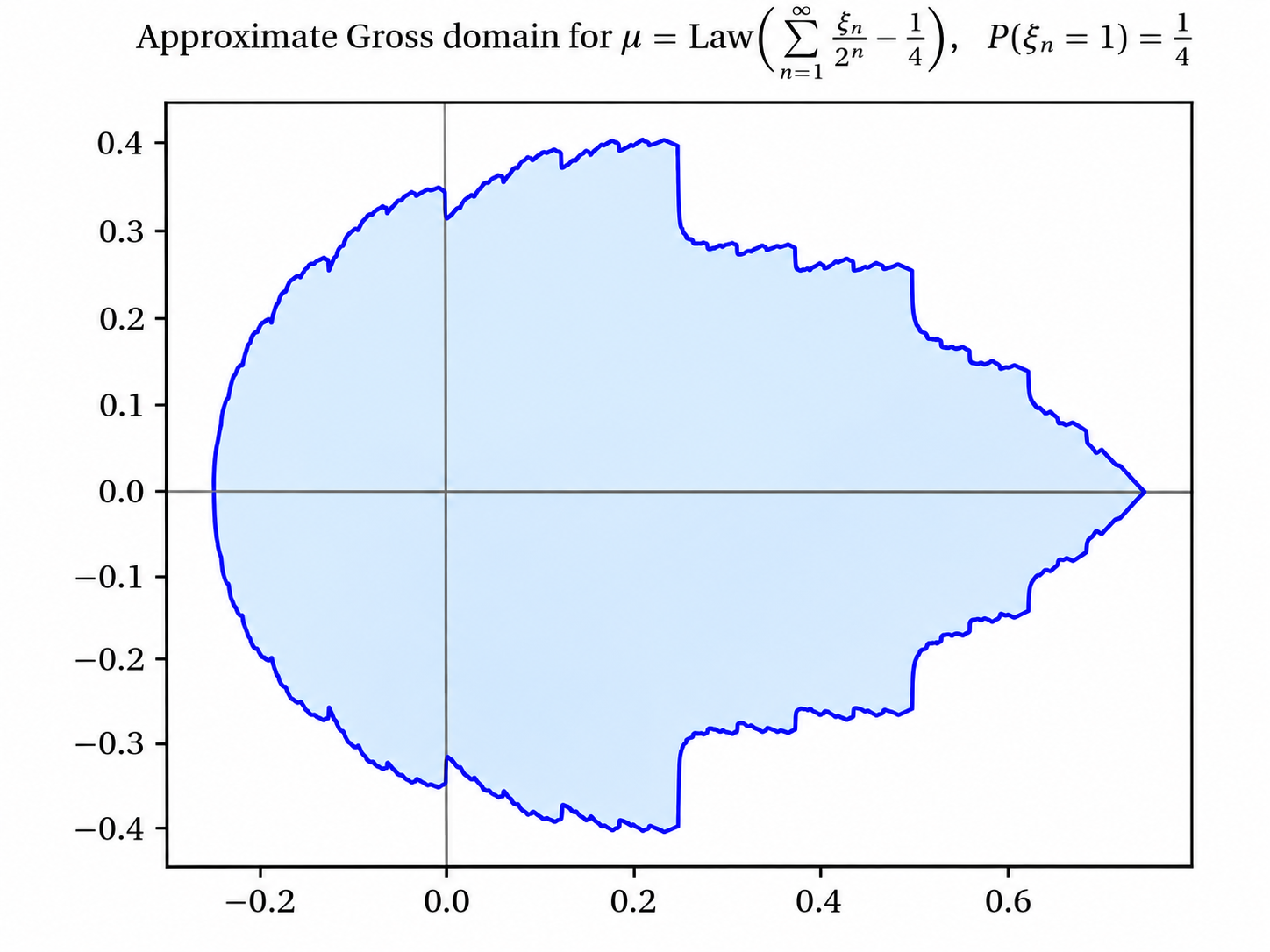} 
\par\end{centering}
\caption{Gross embedding}
\end{figure}

The next result reveals a necessary measure-theoretic condition for
the existence of finite-perimeter embeddings: the prescribed exit
law must dominate Lebesgue measure on its support, although it need
not itself be absolutely continuous.

\begin{theorem}
\label{thm:necessary-density} If the $\operatorname{CPOP}$ is nondegenerate
then 
\[
\mathcal{L}^{1}|_{\operatorname{supp}(\mu)}\ll\mu.
\]
Equivalently, if 
\[
\mu(dx)=\rho(x)\,dx+\mu_{s}(dx)
\]
is the Lebesgue decomposition of $\mu$, then 
\[
\rho(x)>0\quad\text{for Lebesgue-almost every }x\in\operatorname{supp}(\mu).
\]
\end{theorem} 
\begin{proof}
Denote by $S$ the $\operatorname{supp}(\mu)$. Since $m_{c}(\mu)<\infty$,
there exists $U\in\mathcal{C}(\mu)$ with finite conformal perimeter.
Let $f_{U}:\mathbb{D}\to U$ be its normalized conformal map. By definition,
$f_{U}'\in\mathbf{H}^{1}$. It follows from the classical boundary
regularity theorem for Hardy spaces that $f$ extends continuously
to $\overline{\mathbb{D}}$ and that its boundary trace 
\[
\gamma(t)=f_{U}(e^{it}),\qquad0\le t\le2\pi,
\]
is absolutely continuous. In particular, 
\[
u(t)=\operatorname{Re}\gamma(t)
\]
is an absolutely continuous real-valued function. By conformal invariance
of harmonic measure, 
\[
\mu=u_{*}\left(\frac{dt}{2\pi}\right).
\]
Since $u$ is continuous, the support of this pushforward measure
coincides with its image. Hence 
\[
S=u([0,2\pi]),
\]
which is a compact interval. Now let $E\subseteq S$ be a Borel set
such that $\mu(E)=0$. By the preceding identity, 
\[
\mathcal{L}^{1}(u^{-1}(E))=2\pi\mu(E)=0.
\]
Since absolutely continuous functions map Lebesgue-null sets to Lebesgue-null
sets, 
\[
\mathcal{L}^{1}\bigl(u(u^{-1}(E))\bigr)=0.
\]
But $E\subseteq S=u([0,2\pi])$, so 
\[
u(u^{-1}(E))=E.
\]
Consequently, 
\[
\mathcal{L}^{1}(E)=0,
\]
which proves 
\[
\mathcal{L}^{1}|_{S}\ll\mu.
\]
Finally, write $\mu=\rho\,\mathcal{L}^{1}+\mu_{s}$, where $\mu_{s}$
is supported on a Lebesgue-null set $N$. If $\rho$ vanished on a
subset $E\subseteq S$ of positive Lebesgue measure, then 
\[
\mathcal{L}^{1}(E\setminus N)>0,\qquad\mu(E\setminus N)=0,
\]
contradicting $\mathcal{L}^{1}|_{S}\ll\mu$. Therefore $\rho>0$ almost
everywhere on $S$.
\end{proof}

\begin{remark} The preceding condition is not equivalent to absolute
continuity of $\mu$ with respect to Lebesgue measure. Indeed, finite
conformal perimeter does not exclude a singular component: the real-part
exit distribution of Brownian motion from a square has atoms corresponding
to its vertical sides, although the square has finite conformal perimeter.
Thus, the necessary condition concerns the positivity of the absolutely
continuous component of $\mu$, rather than the absence of a singular
component. \end{remark}

Theorem~\ref{thm:necessary-density} immediately implies a Hausdorff-dimension
result. Recall that the Hausdorff dimension of $\mu$ (a probability
measure in general) is defined by 
\[
\dim_{\mathcal{H}}\mu=\inf\left\{ \dim_{\mathcal{H}}(E):\mu(E)=1\right\} .
\]
\begin{corollary} \label{cor:hausdorff-dimension} If $m_{c}(\mu)<\infty$,
then 
\[
\dim_{\mathcal{H}}\mu=1.
\]
In particular, 
\[
\dim_{\mathcal{H}}\mu<1\quad\Longrightarrow\quad m_{c}(\mu)=+\infty.
\]
\end{corollary}
\begin{proof}
By Theorem~\ref{thm:necessary-density}, $S=\operatorname{supp}\mu$
is a nondegenerate compact interval and $\mathcal{L}^{1}|_{S}\ll\mu$.
Thus every Borel set $E$ with $\mu(E)=1$ satisfies 
\[
\mathcal{H}^{1}(E\cap S)=\mathcal{L}^{1}(S)>0.
\]
Hence $\dim_{\mathcal{H}}(E)\geq1$. Since $\mu$ is supported on
$\mathbb{R}$, the conclusion follows. 
\end{proof}

\subsection{Sufficient conditions for existence of a perimeter minimizer}

In \cite{Energy}, the Skorokhod energy of $U\in\mathscr{C}(\mu)$
is defined by 
\[
\Lambda(U)=\frac{1}{4}\Vert f_{U}'\Vert^{2}_{\mathbf{H}^{2}},
\]
and it is shown that 
\[
\Lambda(U^{G}_{\mu})\leq\Lambda(U).
\]
The quantity $\Lambda(U^{G}_{\mu})$ is referred to as the Skorokhod
energy of $\mu$ and denoted by $\Lambda_{\mu}$. When $\mu$ has
a density $\rho$ and no gaps, then $\Lambda_{\mu}$ admits the integral
formula 
\[
\Lambda_{\mu}=\frac{1}{2\pi^{2}}\int_{\operatorname{supp}(\mu)}\frac{dx}{\rho(x)}.
\]
\begin{proposition}[Bounds for the optimal conformal perimeter]
\label{prop:perimeter-bounds} We have 
\begin{equation}
2\pi\sqrt{2\operatorname{Var}(\mu)}\leq m_{c}(\mu)\leq4\pi\sqrt{\Lambda_{\mu}}.\label{eq:perimeter-bounds}
\end{equation}
Moreover, both bounds are sharp. In particular, if $\Lambda_{\mu}$
is finite then \ref{optimization problem} is nondegenerate.\end{proposition} 
\begin{proof}
Let $U\in\mathscr{C}(\mu)$ of finite perimeter, and write 
\[
f_{U}(z)=\sum^{\infty}_{n=1}c_{n}z^{n}.
\]
In particular, $U$ is bounded by \ref{bishopestimate}. Hence $\mathbf{E}(\tau_{U})<\infty$,
and so $f_{U}\in\mathbf{H}^{2}$ by Burkholder equivalence theorem
\cite{Burkholder1977}. In this case
\[
\mathbf{E}(\tau_{U})=\frac{1}{2}\sum^{\infty}_{n=1}|c_{n}|^{2}=\operatorname{Var}(\mu).
\]
See \cite{Boudabra2026}. On the other hand, the area formula for
univalent functions gives 
\[
\operatorname{Area}(U)=\pi\sum^{\infty}_{n=1}n|c_{n}|^{2}\ge2\pi\operatorname{Var}(\mu).
\]
Combining the isoperimetric inequality with the domination of geometric
perimeter by conformal perimeter yields 
\[
\begin{aligned}\operatorname{Per}_{c}(U) & \ge2\sqrt{\pi\operatorname{Area}(U)}\\
 & \ge2\pi\sqrt{2\operatorname{Var}(\mu)}.
\end{aligned}
\]
Taking the infimum over all admissible domains proves the lower bound.
For the upper bound, Baunyakovsky-Cauchy--Schwarz inequality applied
to Gross map $g$ yields 
\[
\begin{aligned}m_{c}(\mu) & \le\operatorname{Per}_{c}(U^{G}_{\mu})\\
 & =2\pi\|g'\|_{\mathbf{H}^{1}}\\
 & \le2\pi\|g'\|_{\mathbf{H}^{2}}\\
 & =4\pi\sqrt{\Lambda_{\mu}}.
\end{aligned}
\]
The sharpness of the bounds follow for example from the centered arcsine
distribution on $[-1,1]$, whose law is 
\[
d\mu_{\mathrm{as}}(x)=\frac{\mathbf{1}_{(-1,1)}(x)}{\pi\sqrt{1-x^{2}}}\,dx.
\]
where in such a case we have 
\[
2\pi\sqrt{2\operatorname{Var}(\mu_{as})}=m_{c}(\mu_{as})=2\pi=4\pi\sqrt{\Lambda_{\mu_{as}}}.
\]
If $\Lambda_{\mu}$ is finite then Gross embedding would have a finite
conformal perimeter, hence the existence theorem \ref{thm:perimeter-existence}
applies to $\mu$. This completes the proof. 
\end{proof}

\begin{corollary} If $\mu$ has a density $\rho$ bounded below by
a positive threshold then \ref{optimization problem} is nondegenerate.
\end{corollary}

We now establish a refined sufficient condition for the existence
of finite-conformal-perimeter embeddings. Based on logarithmic integrability,
this condition is strictly weaker than the one arising from finite
Skorokhod energy, thereby extending the class of distributions for
which the perimeter optimization problem admits a finite minimizer.

\begin{proposition}\label{prop:log-density-perimeter}
Suppose that $\mu$ has no gaps and a density $\rho>0$. If 
\[
\int\log\left(1+\frac{1}{\rho(x)}\right)\,dx<\infty,
\]
then the Gross embedding has finite conformal perimeter. In particular,
$m_{c}(\mu)<\infty$, and the conformal-perimeter infimum is attained.
\end{proposition} 
\begin{proof}
Let $Q$ be the quantile function of $\mu$. Since $\rho>0$ almost
everywhere, $Q$ is absolutely continuous and satisfies 
\[
Q'(t)=\frac{1}{\rho(Q(t))}\qquad\text{for almost every }t\in(0,1).
\]
The symmetric decreasing boundary datum defining the Gross conformal
map is 
\[
q(\theta)=Q\left(1-\frac{|\theta|}{\pi}\right),\qquad-\pi\leq\theta\leq\pi.
\]
In particular, $q$ is absolutely continuous and 
\[
|q'(\theta)|=\frac{1}{\pi}Q'\left(1-\frac{|\theta|}{\pi}\right)
\]
for almost every $\theta$. A change of variables gives 
\[
\begin{aligned} & \int^{\pi}_{-\pi}|q'(\theta)|\log\bigl(1+|q'(\theta)|\bigr)\,d\theta\\
 & \qquad=2\int^{1}_{0}Q'(t)\log\left(1+\frac{Q'(t)}{\pi}\right)\,dt\\
 & \qquad=2\int^{b}_{a}\log\left(1+\frac{1}{\pi\rho(x)}\right)\,dx\\
 & \qquad<\infty.
\end{aligned}
\]
Thus $q'\in L\log L(\mathbb{T})$. By the classical Zygmund inequality
for the periodic Hilbert transform $\mathcal{H}$ \cite{Yee}, 
\[
q'\in L\log L(\mathbb{T})\quad\Longrightarrow\quad\mathcal{H}q'\in L^{1}(\mathbb{T}).
\]
The Gross conformal map has boundary values 
\[
g^{*}=q+i\mathcal{H}q.
\]
Since the Hilbert transform commutes with distributional differentiation,
\[
\frac{d}{d\theta}g^{*}=q'+i\mathcal{H}q'\in L^{1}(\mathbb{T}).
\]
Writing $P_{r}$ for the Poisson kernel, we obtain 
\[
ire^{i\theta}g'(re^{i\theta})=P_{r}*\bigl(q'+i\mathcal{H}q'\bigr)(\theta).
\]
Consequently, 
\[
r\int^{\pi}_{-\pi}|g'(re^{i\theta})|\,d\theta\leq\int^{\pi}_{-\pi}|q'+i\mathcal{H}q'|\,d\theta<\infty.
\]
Letting $r\uparrow1$ and using the monotonicity of Hardy integral
means yields 
\[
g'\in\mathbf{H}^{1}.
\]
Therefore, 
\[
\operatorname{Per}_{c}(U^{G}_{\mu})=2\pi\|g'\|_{\mathbf{H}^{1}}<\infty.
\]
Since $U^{G}_{\mu}\in\mathcal{C}(\mu)$, we conclude that $m_{c}(\mu)<\infty$.
Attainment follows from the conformal-perimeter existence theorem
\ref{thm:perimeter-existence}. 
\end{proof}

The logarithmic condition accommodates a remarkably broad class of
distributions, including densities with strong degeneracies. For instance,
the centered probability densities 
\[
\rho_{\alpha}(x)=\frac{\alpha+1}{2}|x|^{\alpha}\mathbf{1}_{(-1,1)}(x),\qquad\alpha>0,
\]
satisfy the condition for every $\alpha>0$, although the corresponding
Skorokhod-energy integrability condition fails whenever $\alpha\geq1$.
More strikingly, consider 
\[
\rho_{\alpha}(x)=c_{\alpha}\exp\left(-\frac{1}{|x|^{\alpha}}\right)\mathbf{1}_{(-1,1)}(x),\qquad0<\alpha<1,
\]
where $c_{\alpha}$ is the normalizing constant and $\rho_{\alpha}(0)=0$.
Since 
\[
\log\left(1+\frac{1}{\rho_{\alpha}(x)}\right)\sim\frac{1}{|x|^{\alpha}}\qquad(x\to0),
\]
the logarithmic condition is satisfied. Nevertheless, 
\[
\int^{1}_{-1}\rho_{\alpha}(x)^{-\varepsilon}\,dx=+\infty\qquad\text{for every }\varepsilon>0.
\]
Thus, the logarithmic criterion guarantees the existence of a finite-conformal-perimeter
minimizer even for distributions whose densities vanish faster than
any polynomial and for which every positive inverse-power integrability
condition fails.

\section{A family of distributions where Gross embedding is not optimal }\label{A-family-of}

Let $U$ be an embedding in $\mathscr{C}(\mu)$. For $s\in(0,1)$
the fractional boundary energy $\mathscr{E}_{s}(U)$ of $U$ is defined
by 
\[
\mathscr{E}_{s}(U)=\sum^{+\infty}_{n=1}n^{2s}|c_{n}|^{2}
\]
where the $c_{n}$'s are the coefficients of $f_{U}$. In \cite{Boudabra2026,Pruss},
It is proved that Gross embedding represents the unique minimizer
of $\mathscr{E}_{s}$. In particular, Gross embedding has the least
area which corresponds to $\pi\mathscr{E}_{1/2}$. This further reinforced
our earlier conjecture that the Gross embedding would also minimize
perimeter. However, we were unable to show that. Recently, we thought
to investigate the other scenario: what if Gross embedding does not
the minimal perimeter? So we starting by selecting random univalent
functions $f$ and compute the perimeter of their ranges $f(\mathbb{D})$,
and then compare it with their Brownian symmetrizations. Surprisingly,
some choices show that 
\[
\Per_{c}(f(\mathbb{D}))<\Per_{c}(\mathfrak{B}(f(\mathbb{D})))!
\]
Repeating the computations yielded the same comparison. We succeeded
to turn the observations into a rigorous work. The idea is to use
a perturbation approach. Before constructing the \textquotedbl counterexamples\textquotedbl ,
we begin by some technical tools.

The classical Hardy point-evaluation inequality is 
\begin{equation}
|h(w)|\le\frac{\|h\|_{\mathbf{H}^{1}}}{1-|w|^{2}},\qquad h\in\mathbf{H}^{1},\quad w\in\D.\label{eq:classical-hardy}
\end{equation}
Recall its elementary mechanism. For $H(z)=\sum_{n\ge0}a_{n}z^{n}\in\mathbf{H}^{2}$,
Parseval and Cauchy--Schwarz give 
\[
|H(w)|^{2}\le\left(\sum_{n\ge0}|a_{n}|^{2}\right)\left(\sum_{n\ge0}|w|^{2n}\right)=\frac{\|H\|^{2}_{\mathbf{H}^{2}}}{1-|w|^{2}}.
\]
If $h$ is zero-free, choose an analytic square root $\mathbf{H}^{2}=h$;
then $\|H\|^{2}_{\mathbf{H}^{2}}=\|h\|_{\mathbf{H}^{1}}$, proving
\eqref{eq:classical-hardy}. For general $h$, the classical Riesz
factorization $h=AB$ with $A,B\in\mathbf{H}^{2}$ and $\|A\|_{\mathbf{H}^{2}}\|B\|_{\mathbf{H}^{2}}=\|h\|_{\mathbf{H}^{1}}$
gives the same estimate by applying the $\mathbf{H}^{2}$ inequality
twice; see~\cite{Duren1970}. It turns out that if $h$ is even then
we can refine the Hardy estimation \ref{eq:classical-hardy} a tiny
bit more.

\begin{lemma}$\dagger$ \label{lem:hardy} If $h\in\mathbf{H}^{1}$
is even, then, for $0<r<1$, 
\begin{equation}
\|h\|_{\mathbf{H}^{1}}\ge(1-r^{4})|h(r)|.\label{eq:even-hardybnd}
\end{equation}
\end{lemma}

\begin{corollary}If $h$ is odd univalent function on the unit disc
then 
\begin{equation}
\Per_{c}(h(\D))\ge2\pi(1-r^{4})|h'(r)|.\label{eq:perimeter-point}
\end{equation}

\end{corollary}

\begin{lemma} $\dagger$ \label{lem:boundary} Suppose $q$ is a
real-valued, even, $2\pi$-periodic $\mathcal{C}^{3}$ function of
zero mean, strictly decreasing on $[0,\pi]$. Its analytic completion
\eqref{eq:gross-series} extends to a $C^{1}$ map on $\overline{\D}$,
is injective there, and maps $\D$ onto a bounded rectifiable Jordan
domain. In particular, \ref{eq:perimeters-equal} applies to $g$.
\end{lemma}

\subsection{A cubic family and the perimeter comparison}

\label{sec:family}

In this section we will provide examples where $\mathfrak{B}(U)$
has a larger perimeter than $U$. The approach is to perturb the arcsine
distribution. Consider 
\begin{equation}
F_{b}(z)=z+\frac{b}{3}z^{3},\qquad W_{a,b}=F_{b}(\D)-F_{b}(ia),\qquad0<a<1,\quad0<b\le1.\label{eq:family}
\end{equation}
Write $G_{a,b}=\Bsym(W_{a,b})$ and 
\begin{equation}
k=\frac{1-a^{2}}{1+a^{2}}\in(0,1).\label{eq:k}
\end{equation}

\begin{theorem} $\dagger$ \label{thm:family} The domains $W_{a,b}$
and $G_{a,b}$ are bounded rectifiable Jordan domains containing $0$
and having the same symmetric real exit law. Their perimeters satisfy
\begin{align}
\Per_{c}(W_{a,b}) & \le2\pi(1-b)+8b,\label{eq:family-upper}\\
\Per_{c}(G_{a,b}) & >\frac{8}{15}(1+a^{2})(5+3b).\label{eq:family-lower}
\end{align}
In particular, whenever 
\begin{equation}
\frac{8}{15}(1+a^{2})(5+3b)>2\pi(1-b)+8b,\label{condition ab}
\end{equation}
one has 
\begin{equation}
\Per_{c}(G_{a,b})>\Per_{c}(W_{a,b}).\label{reverse inequa per}
\end{equation}
For $b=1$, the perimeter in \eqref{eq:family-upper} is exactly $8$.
For $0<b<1$, both boundaries are real-analytic Jordan curves. \end{theorem}

\begin{remark} Note that Theorem \ref{thm:family} describes a family
of counterexamples; the prescribed exit law varies with the parameters.
Nevertheless, it does not compute $m_{c}(\mu_{a,b})$ where $\mu_{a,b}$
is the law of $\operatorname{Re}(Z_{\tau_{W_{a,b}}})$.

\end{remark}

\begin{figure}[ht]
\centering \begin{tikzpicture}[scale=3, >=stealth]

\pgfmathsetmacro{\bz}{(30*pi-80)/(30*pi-72)}
\pgfmathsetmacro{\aone}{sqrt(7/8)}

\fill[blue!20]
  plot[domain=\bz:1, samples=100, variable=\t]
  ({sqrt(15*(2*pi+(8-2*pi)*\t)/(8*(5+3*\t))-1)}, {\t})
  -- (1,1) -- cycle;

\draw[->] (0,0) -- (1.15,0) node[right] {$a$};
\draw[->] (0,0) -- (0,1.15) node[above] {$b$};

\draw[dashed] (1,0) -- (1,1);
\draw[thick]  (0,1) -- (1,1);

\draw[thick, red!70!black, dashed, domain=\bz:1, samples=100, variable=\t]
  plot ({sqrt(15*(2*pi+(8-2*pi)*\t)/(8*(5+3*\t))-1)}, {\t});

\draw (1,0) -- (1,-0.015) node[below] {$1$};
\draw (0,1) -- (-0.015,1) node[left] {$1$};
\draw (0,\bz) -- (-0.015,\bz) node[left] {$b_0$};
\draw[dotted] (0,\bz) -- (1,\bz);
\draw (\aone,1) -- (\aone,1.015) node[above, font=\small] {$\sqrt{7/8}$};

\node at (0.93,0.83) {$\mathcal{R}$};

\end{tikzpicture} \caption{The region $\mathcal{R}=\{(a,b):\sqrt{7/8}<a<1,\ f(a)<b\le1\}$ of
parameters satisfying $\frac{8}{15}(1+a^{2})(5+3b)>2\pi(1-b)+8b$,
where $f(a)=\frac{30\pi-40-40a^{2}}{30\pi-96+24a^{2}}$. In the region
$\mathcal{R}$, $\Per(G_{a,b})>\Per(W_{a,b})$ is guaranteed. The
dashed curve $b=f(a)$ and the line $a=1$ are excluded; the segment
$b=1$ is included. The curve meets $b=1$ at $a=\sqrt{7/8}\approx0.935$
and $a=1$ at $b_{0}=\frac{30\pi-80}{30\pi-72}\approx0.640$.}
\label{fig:region-perimeter} 
\end{figure}
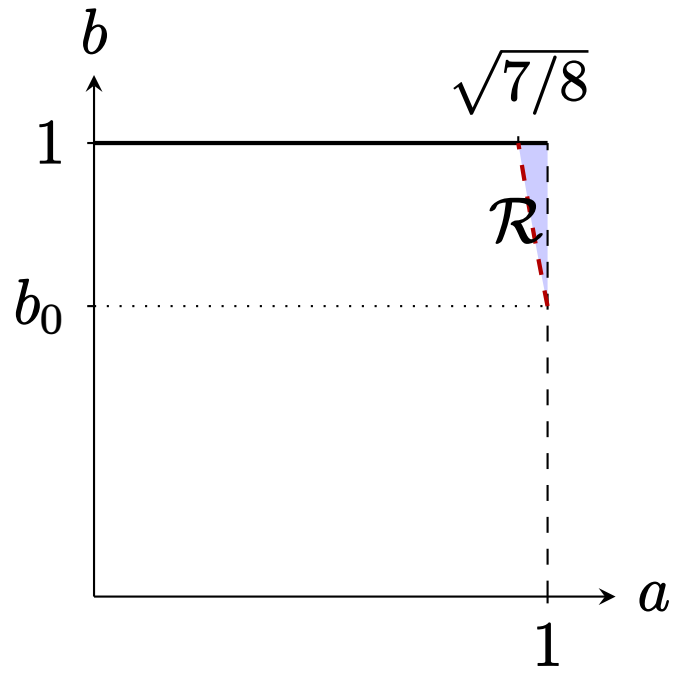

For illustration, we give two worked examples. For the first example
take $b=1$ and $a=2\sqrt{2}/3$. Then $a^{2}=8/9$ and 
\[
F_{1}(ia)=ia\left(1-\frac{a^{2}}{3}\right)=\frac{38\sqrt{2}}{81}\,i.
\]
Thus $W_{a,1}=U$, and Proposition~\ref{thm:family} gives 
\[
\Per(U)=8,\qquad\Per(\Bsym(U))>\frac{8}{15}\frac{17}{9}\,8=\frac{1088}{135}.
\]
In this case $k=1/17$, and the quantile in \eqref{eq:family-quantile}
reduces to 
\begin{equation}
q(t)=\frac{4\cos^{3}t}{3(\cos^{2}t+289\sin^{2}t)^{3/2}}.\label{eq:cubic-quantile}
\end{equation}

For the second example take $b=3/4$ and $a=2\sqrt{6}/5$. Then $a^{2}=24/25$,
$k=1/49$, and 
\[
F_{3/4}(ia)=ia\left(1-\frac{a^{2}}{4}\right)=\frac{38\sqrt{6}}{125}\,i.
\]
Hence $W_{a,3/4}=V$, and both boundaries are real-analytic. The bounds
in Proposition~\ref{thm:family} become 
\[
\Per(V)\le6+\frac{\pi}{2}<\frac{53}{7},\qquad\Per(\Bsym(V))>\frac{8}{15}\frac{49}{25}\frac{29}{4}=\frac{2842}{375}.
\]
Here $\pi<22/7$ is sufficient; it follows, for instance, from the
positive integral 
\[
\frac{22}{7}-\pi=\int^{1}_{0}\frac{x^{4}(1-x)^{4}}{1+x^{2}}\,dx>0.
\]
Finally, 
\[
\frac{2842}{375}-\frac{53}{7}=\frac{19}{2625}>0.
\]
The supports and symmetry of the two exit laws follow from \eqref{eq:family-quantile}.

\subsection{Normalization and exact quadratic norms}

The conformal perimeter is independent of the choice of conformal
parametrization. Derivative $H^{p}$ norms for $p\ne1$, however,
need not be invariant under precomposition by arbitrary disc automorphisms.
We therefore normalize at the Brownian starting point. For the family
\eqref{eq:family}, put 
\begin{equation}
\phi_{a}(z)=\frac{z+ia}{1-iaz},\qquad f_{a,b}(z)=F_{b}(\phi_{a}(z))-F_{b}(ia)\label{eq:normalized-map}
\end{equation}
where again $F_{b}(z)=z+bz^{3}/3$ (the same function in Section \ref{sec:family}).
Then $f_{a,b}(0)=0$ and $f_{a,b}'(0)=(1-ba^{2})(1-a^{2})>0$. The
Gross map $g_{a,b}$ also fixes $0$ and has positive derivative there.
All norms compared in this section use these maps.\\
 \\

\textbf{Comment 1. Persistence near $p=1$ and an equality exponent}.

Let $a,b$ for which \ref{condition ab} is satisfied. Then 
\begin{equation}
\|f_{a,b}'\|_{\mathbf{H}^{1}}<\|g_{a,b}'\|_{\mathbf{H}^{1}}.\label{eq:p=00003D00003D1}
\end{equation}
By the continuity of the $\mathbf{H}^{p}$ norms in $p>0$, the inequality
\ref{eq:p=00003D00003D1} persists near $p=1^{+}$. There exists $\eta>0$
such that 
\begin{equation}
\|f_{a,b}'\|_{\mathbf{H}^{p}}<\|g_{a,b}'\|_{\mathbf{H}^{p}}\qquad(1\le p<1+\eta).\label{eq:hp-failure}
\end{equation}
On the other hand, as the Skorokhod energy is exclusively minimal
for Gross embedding, one has 
\[
\|g_{a,b}'\|_{\mathbf{H}^{2}}<\|h_{a,b}'\|_{\mathbf{H}^{2}}.
\]
For that same pair, there also exists at least one $p_{*}\in(1,2)$
such that 
\begin{equation}
\|g'\|_{\mathbf{H}^{p_{*}}}=\|f'\|_{\mathbf{H}^{p_{*}}}.\label{eq:hp-equality}
\end{equation}
In particular, derivative $\mathbf{H}^{p}$ contraction under Gross
symmetrization fails on an interval immediately to the right of $1$,
even for domains with real-analytic boundaries. This statement neither
identifies $p_{*}$ nor proves its uniqueness. It also does not assert
a universal transition exponent independent of the exit law.\\
 \\

\textbf{Comment 2. A failed rearrangement inequality}

Let $u(t)=\operatorname{Re}f(e^{it})$ for either of the normalized
maps above, and let $u^{*}$ be its even decreasing rearrangement.
Then $u^{*}=q$, the Gross boundary datum. Define the periodic Hilbert
transform on Fourier modes by 
\[
\cH(\cos nt)=\sin nt,\qquad\cH(\sin nt)=-\cos nt,\qquad\cH1=0.
\]
Since $f(0)=0$, analytic completion gives $f(e^{it})=u(t)+i\cH u(t)$.
Differentiation and the arc-length formula show that 
\begin{equation}
\mathcal{J}(u):=\int^{2\pi}_{0}\sqrt{u'(t)^{2}+(\cH u'(t))^{2}}\,dt=2\pi\|f'\|_{\mathbf{H}^{1}}.\label{eq:hilbert-functional}
\end{equation}
The same identity holds for $u^{*}$ and $g$. Theorem \ref{thm:family}
therefore gives the strict reverse comparison 
\begin{equation}
\mathcal{J}(u^{*})>\mathcal{J}(u).\label{eq:failed-rearrangement}
\end{equation}
This disproves a general rearrangement-decreasing inequality for $\mathcal{J}$,
even when both analytic completions are univalent and, using $V$,
have regular analytic boundary traces.

In contrast, the scalar total variations agree. The function $r_{b}$
in \eqref{eq:real-boundary} traverses $[-1-b/3,1+b/3]$ monotonically
on each semicircle; the boundary automorphism in \eqref{eq:normalized-map}
does not change its total variation. The same is true of $q$. Hence
\[
\int^{2\pi}_{0}|u'(t)|\,dt=\int^{2\pi}_{0}|(u^{*})'(t)|\,dt=4\left(1+\frac{b}{3}\right).
\]
The failure in \eqref{eq:failed-rearrangement} concerns the full
complex boundary, not an increase in variation of its real coordinate.\\
 \\

\textbf{Comment 3. The isoperimetric ratio can increase.}

For a bounded rectifiable Jordan domain $W$, define 
\[
I(W)=\frac{\Per(W)^{2}}{4\pi\Area(W)}.
\]
The area-nonincreasing property 
\begin{equation}
\Area(\Bsym(W))\le\Area(W)\label{eq:area-monotone}
\end{equation}
is the case $s=1/2$ of~\cite[Theorem~3.4]{Boudabra2026}, since
\[
\Area(W)=\pi\sum_{n\ge1}n|a_{n}|^{2}
\]
for the normalized univalent map $f_{W}(z)=\sum_{n\ge1}a_{n}z^{n}$.
This cited area result is an additional input to the next consequence;
it is not used in the counterexamples or Hardy-norm calculations.\\
 \\

\begin{corollary}[Isoperimetric amplification in the cubic family]
$\dagger$ \label{cor:isoperimetric} For $0<a<1$ and $0<b\leq1$,
define 
\[
R(a,b)=\frac{\frac{8}{15}(1+a^{2})(5+3b)}{2\pi(1-b)+8b}.
\]
Then the domains in \eqref{eq:family} satisfy 
\begin{equation}
\frac{I(G_{a,b})}{I(W_{a,b})}>R(a,b)^{2}.\label{eq:iso-family}
\end{equation}
In particular, 
\begin{equation}
2\pi(1-b)+8b\leq\frac{8}{15}(1+a^{2})(5+3b)\quad\Longrightarrow\quad I(G_{a,b})>I(W_{a,b}).\label{eq:iso-family-condition}
\end{equation}

\end{corollary}

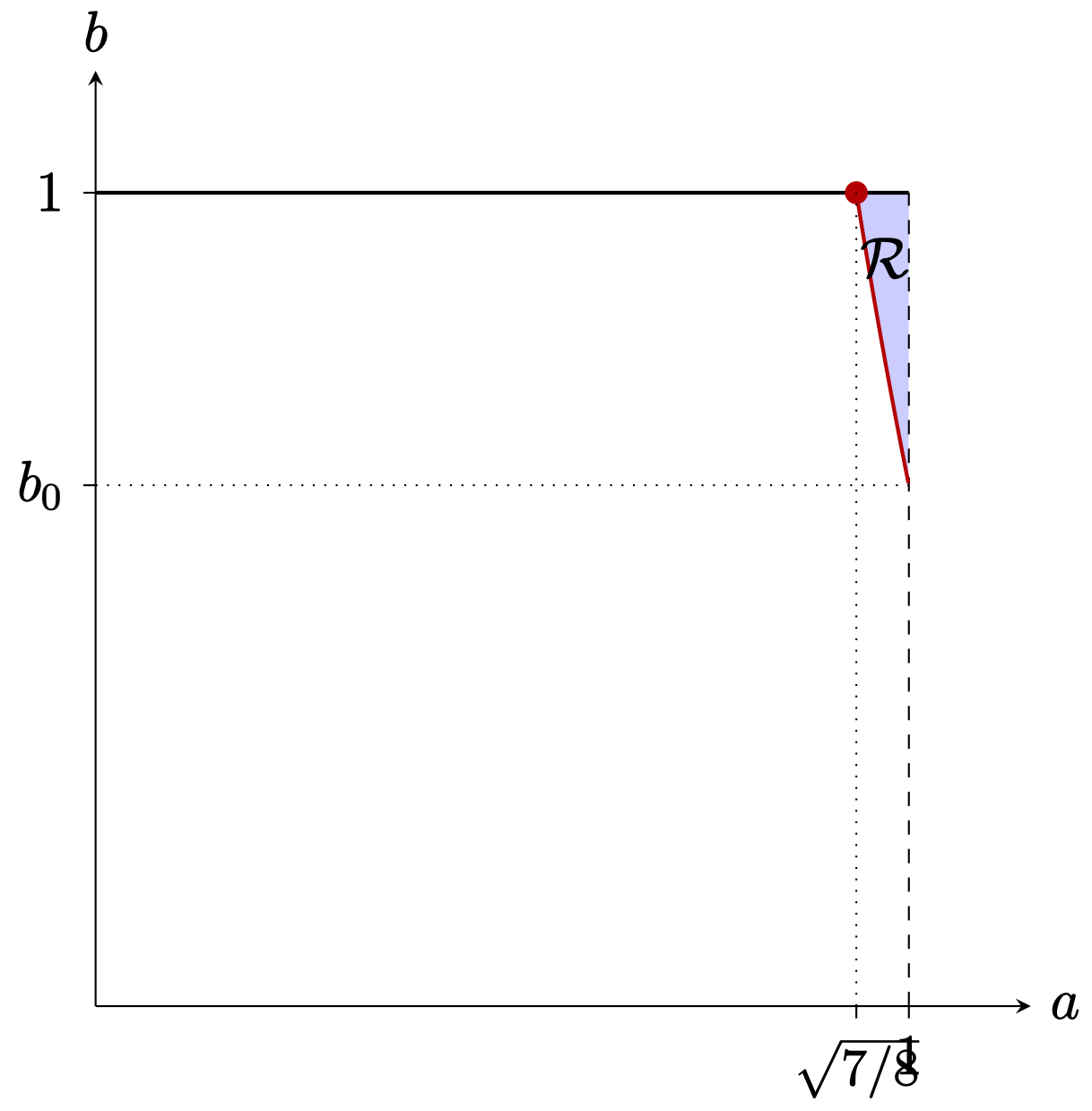
\begin{figure}[htbp]
\centering \begin{tikzpicture}[scale=6, >=stealth]

\pgfmathsetmacro{\bz}{(30*pi-80)/(30*pi-72)}   
\pgfmathsetmacro{\aone}{sqrt(7/8)}             

\fill[blue!20]
  plot[domain=\aone:1, samples=100, variable=\t]
  ({\t}, {(30*pi-40-40*\t*\t)/(30*pi-96+24*\t*\t)})
  -- (1,1) -- cycle;

\draw[->] (0,0) -- (1.15,0) node[right] {$a$};
\draw[->] (0,0) -- (0,1.15) node[above] {$b$};

\draw[dashed] (1,0) -- (1,1);
\draw[thick]  (0,1) -- (1,1);

\draw[thick, red!70!black, domain=\aone:1, samples=100, variable=\t]
  plot ({\t}, {(30*pi-40-40*\t*\t)/(30*pi-96+24*\t*\t)});

\fill[red!70!black] (\aone,1) circle (0.4pt);

\draw (1,0) -- (1,-0.015) node[below] {$1$};
\draw (0,1) -- (-0.015,1) node[left] {$1$};
\draw (0,\bz) -- (-0.015,\bz) node[left] {$b_0$};
\draw[dotted] (0,\bz) -- (1,\bz);
\draw (\aone,0) -- (\aone,-0.015) node[below, font=\small] {$\sqrt{7/8}$};
\draw[dotted] (\aone,0) -- (\aone,1);

\node at (0.97,0.92) {$\mathcal{R}$};

\end{tikzpicture} \caption{The region $\mathcal{R}=\{(a,b):\sqrt{7/8}\le a<1,\ f(a)\le b\le1\}$
of parameters satisfying $\Big(\frac{\frac{8}{15}(1+a^{2})(5+3b)}{2\pi(1-b)+8b}\Big)^{2}\ge1$,
where $f(a)=\frac{30\pi-40-40a^{2}}{30\pi-96+24a^{2}}$. The curve
$b=f(a)$ and the segment $b=1$ are included; the line $a=1$ is
excluded. The curve meets $b=1$ at $a=\sqrt{7/8}\approx0.935$ and
$a=1$ at $b_{0}=\frac{30\pi-80}{30\pi-72}\approx0.640$.}
\label{fig:region-R} 
\end{figure}

For the domain $U=W_{\frac{2\sqrt{2}}{3},1}$ , this gives $R(\frac{2\sqrt{2}}{3},1)=136/135$
and hence 
\begin{equation}
\frac{I(\Bsym(U))}{I(U)}>1.\label{eq:iso-gap}
\end{equation}

\begin{remark} Condition \eqref{eq:iso-family-condition} is sufficient,
not necessary, for an increase of the isoperimetric ratio. When it
fails, estimate \eqref{eq:iso-family} alone does not determine whether
$I(G_{a,b})>I(W_{a,b})$. Thus Brownian symmetrization can increase
perimeter and the isoperimetric ratio while not increasing area. This
does not conflict with the classical isoperimetric inequality; it
shows that this exit-law-preserving operation need not improve its
ratio.

\end{remark}

\section{Case where Gross embedding is the minimal perimeter }\label{Case-where-Gross}

The preceding counterexamples show that Gross symmetrization need
not decrease perimeter. Nevertheless, we now give a sufficient geometric
condition for perimeter contraction: convexity of the symmetrized
domain. Unlike comparisons restricted to a prescribed family of boundary
parametrizations, the result below applies to every bounded rectifiable
Jordan competitor.

\subsection{The general perimeter comparison}

The convex hull of a domain $U$, denoted by $\operatorname{conv}(U)$
is the smallest convex region containing $U$. A popular result from
convex geometry is that the perimeter of $\operatorname{conv}(U)$
does not exceed that of $U$ \cite{Focardi}. Recall that the convex
hull of an open planar domain $U$ is itself open. In particular,
\[
\operatorname{conv}(U)=\operatorname{int}(\operatorname{conv}(\overline{U})).
\]
Moreover, every bounded convex domain has a rectifiable Jordan boundary,
and hence its conformal perimeter coincides with its geometric perimeter:
\[
\operatorname{Per}_{c}(V)=\operatorname{Per}_{g}(V)
\]
for every bounded convex domain $V$. Consequently, the classical
perimeter inequalities for convex bodies apply equally to their conformal
perimeters which are the choice for this manuscript.

\begin{theorem}[Convex-hull lower bound and uniqueness] \label{thm:convex-hull-lower-bound}
Assume Gross embedding $U^{G}_{\mu}$ is bounded, and set 
\[
K_{\mu}=\operatorname{conv}(U^{G}_{\mu}).
\]
Then 
\[
\boxed{\operatorname{Per}_{c}(K_{\mu})\le m_{c}(\mu).}
\]
In particular, if $U^{G}_{\mu}$ is convex, then 
\[
\boxed{m_{c}(\mu)=\operatorname{Per}_{c}(U^{G}_{\mu}).}
\]
Moreover, in this case, $U^{G}_{\mu}$ is the unique minimizer: for
every $U\in\mathcal{C}(\mu)$, 
\[
\operatorname{Per}_{c}(U)=\operatorname{Per}_{c}(U^{G}_{\mu})\quad\Longleftrightarrow\quad U=U^{G}_{\mu}.
\]
\end{theorem} 
\begin{proof}
If $m_{c}(\mu)=\infty$, the lower bound is immediate. Otherwise,
consider $U\in\mathcal{C}(\mu)$ such that $\operatorname{Per}_{c}(U)<\infty$,
and put 
\[
K_{U}=\operatorname{conv}(U).
\]
Since the normalized conformal map $f_{U}:\mathbb{D}\to U$ satisfies
$f_{U}'\in\mathbf{H}^{1}$, it admits a rectifiable boundary parametrization
of length $\operatorname{Per}_{c}(U)$. In particular, $U$ is bounded.
By the classical convex-hull perimeter inequality, we obtain and 
\[
\operatorname{Per}_{c}(K_{U})\le\operatorname{Per}_{c}(U).
\]
By Pruss's containment theorem \ref{Pruss's cont}, one has 
\[
U^{G}_{\mu}\subseteq\mathscr{S}(U),
\]
Since Steiner symmetrization preserves inclusion, we get 
\[
U^{G}_{\mu}\subseteq\mathscr{S}(U)\subseteq\mathscr{S}(K_{U}).
\]
The set $\mathscr{S}(K_{U})$ is convex, and therefore 
\[
K_{\mu}\subseteq\mathscr{S}(K_{U}).
\]
Using the monotonicity of perimeter under inclusion of planar convex
bodies, together with the classical Steiner perimeter inequality,
we obtain 

\[
\begin{aligned}\operatorname{Per}_{c}(K_{\mu}) & \leq\operatorname{Per}_{c}(\mathscr{S}(K_{U}))\\
 & \leq\operatorname{Per}_{c}(K_{U})\\
 & \leq\operatorname{Per}_{c}(U).
\end{aligned}
\]
Taking the infimum over all admissible domains yields 
\[
\operatorname{Per}_{c}(K_{\mu})\le m_{c}(\mu).
\]
Suppose now that $U^{G}_{\mu}$ is convex. Then $K_{\mu}=U^{G}_{\mu}$,
and since $U^{G}_{\mu}$ is a bounded convex domain, then 
\[
\operatorname{Per}_{c}(K_{\mu})=\operatorname{Per}_{c}(U^{G}_{\mu})\le m_{c}(\mu)\le\operatorname{Per}_{c}(U^{G}_{\mu}),
\]
which proves 
\[
m_{c}(\mu)=\operatorname{Per}_{c}(U^{G}_{\mu}).
\]
It remains to prove uniqueness. Suppose that $U\in\mathcal{C}(\mu)$
is another minimizer. Then 
\[
\operatorname{Per}_{c}(U)=\operatorname{Per}_{c}(U^{G}_{\mu}).
\]
The preceding chain of inequalities becomes 

\[
\operatorname{Per}_{c}(U^{G}_{\mu})=\operatorname{Per}_{c}(\mathscr{S}(K_{U}))=\operatorname{Per}_{c}(K_{U})=\operatorname{Per}_{c}(U).
\]
Since
\[
\overline{U^{G}_{\mu}}\subseteq\overline{\mathscr{S}(K_{U})}
\]
and both sets are compact convex bodies, strict monotonicity of perimeter
under proper inclusion implies
\[
\overline{U^{G}_{\mu}}=\overline{\mathscr{S}(K_{U})}.
\]
Taking interiors, we conclude that
\[
U^{G}_{\mu}=\mathscr{S}(K_{U}).
\]
Since Steiner symmetrization preserves area, we obtain

\[
\operatorname{Area}(U^{G}_{\mu})=\operatorname{Area}(\mathscr{S}(K_{U}))=\operatorname{Area}(K_{U}).
\]
On the other hand, the area inequality for Brownian symmetrization
gives

\[
\operatorname{Area}(U^{G}_{\mu})\leq\operatorname{Area}(U)\leq\operatorname{Area}(K_{U}).
\]
Consequently,

\[
\operatorname{Area}(U^{G}_{\mu})=\operatorname{Area}(U).
\]
Since both domains are images of univalent conformal maps, equality
of their areas is equivalent to equality in the Dirichlet-norm inequality
under Baernstein rearrangement. By Pruss's rigidity theorem \cite{Pruss},
the corresponding conformal maps differ only by a rotation of the
disk parameter. Their images therefore coincide:
\[
U=U^{G}_{\mu}.
\]
This completes the proof of uniqueness.
\end{proof}

A strighforward corollary of Theorem \ref{thm:convex-hull-lower-bound}
is that if the minimizer of \ref{optimization problem} is not Gross
embedding then this latter is necessarily nonconvex.

\subsection{A sufficient condition ensuring the minimality of the perimeter of
Gross embedding.}

There is a classical coefficient criterion for convex analytic functions
\cite{pascu}. If $f=\sum_{n\geq1}a_{n}z^{n}$ is univalent in the
unit disc then a sufficient condition of the convexity of $f(\mathbb{D})$
is the following 
\begin{equation}
\sum_{n\geq2}n^{2}|a_{n}|\leq|a_{1}|.\label{convex criterion}
\end{equation}
Recall the Fourier expansion of the Gross map generated from the distribution
$\mu$: 
\begin{equation}
g_{\mu}(z)=\sum_{n\geq1}a_{n}z^{n},\,\,\,a_{n}=2\int^{1}_{0}q_{\mu}(t)\cos(n\pi t)\,dt.\label{eq:quantile-convexity-coefficients}
\end{equation}
Let 
\[
\xi\sim\operatorname{Kumaraswamy}\left(\frac{1}{2},\frac{1}{2}\right).
\]
Its density on $(0,1)$ is 
\[
\boxed{\rho_{X}(x)=\frac{1}{4\sqrt{x}\sqrt{1-\sqrt{x}}}.}
\]
In particular $\mathbf{E}\xi=8/15$. See \cite{Kumara}. Let $\mu$
be the law of $\xi-8/15$. Computation of the underlying Fourier coefficients
yields 
\[
a_{n}=\text{\ensuremath{\frac{48(-1)^{n+1}}{\pi^{4}n^{4}}}}.
\]
In particular 
\[
\frac{\sum^{\infty}_{n=2}n^{2}|a_{n}|}{|a_{1}|}=\frac{\pi^{2}}{6}-1\approx0.644934<1.
\]
Hence the Gross embedding is convex and it is the minimizer for the
perimeter optimization problem of $\mu$.

\section{Illustrations }\label{Illustrations}

In this section, we provide worked examples associated with several
distributions.

\subsection{The uniform distribution on $(-1,1)$.}

For the uniform probability measure on $(-1,1)$, the Gross embedding
attains the bottom of perimeters in $\mathscr{C}(\mu)$. The increasing
quantile is 
\[
Q_{\mu}(t)=2t-1.
\]
In the increasing-quantile convention, the Gross map is given in \cite{Gross2019}
by 
\[
g(z)=-\frac{8}{\pi^{2}}\sum_{\substack{n\geq1\\
n\ {\rm odd}
}
}\frac{z^{n}}{n^{2}}.
\]
The series converges uniformly on $\overline{\D}$. Its two boundary
branches can be written as 
\[
x\longmapsto x\pm ih(x),\qquad-1\leq x\leq1,
\]
where 
\[
h(x)=\frac{8}{\pi^{2}}\sum_{\substack{n\geq1\\
n\ {\rm odd}
}
}\frac{\sin(n\theta)}{n^{2}},\qquad\theta=\frac{\pi(x+1)}{2}.
\]
The function $h$ is continuous on $[-1,1]$, with 
\[
h(-1)=h(1)=0.
\]
The identity 
\[
\sum_{\substack{n\geq1\\
n\ {\rm odd}
}
}\frac{z^{n}}{n}=\frac{1}{2}\log\frac{1+z}{1-z},\qquad|z|<1,
\]
gives, by radial limits, 
\[
\sum_{\substack{n\geq1\\
n\ {\rm odd}
}
}\frac{\cos(n\theta)}{n}=\frac{1}{2}\log\cot\frac{\theta}{2},\qquad0<\theta<\pi.
\]
Consequently, 
\[
h'(x)=\frac{2}{\pi}\log\cot\frac{\pi(x+1)}{4},\qquad h''(x)=-\frac{1}{\cos(\pi x/2)}<0,\qquad-1<x<1.
\]
Thus $h$ is strictly concave and positive on $(-1,1)$. It follows
that 
\[
U^{G}_{\mu}=\{x+iy:-1<x<1,\ |y|<h(x)\}
\]
is bounded and convex. Theorem \ref{thm:convex-hull-lower-bound}
applies. In particular, the minimum perimeter in the stated class
is 
\[
\Per(U^{G}_{\mu})=2\int^{1}_{-1}\sqrt{1+\frac{4}{\pi^{2}}\left(\log\cot\frac{\pi(x+1)}{4}\right)^{2}}\,dx.
\]
The logarithmic endpoint singularities are integrable, so this quantity
is finite. Using \ref{eq:perimeter-bounds}, one has 
\[
2\pi\sqrt{\frac{2}{3}}\approx5.130\leq\Per_{c}(U^{G}_{\mu})=m_{c}(\mu)\leq4\sqrt{2}\approx5.657
\]
Numerically, we have 
\[
\Per(U^{G}_{\mu})\approx5.3261068066989535.
\]
Note that Gross map satisfies 
\[
\sum_{n\geq1}n^{2}|a_{n}|=\infty.
\]
Nevertheless, Gross embedding is convex as shown. Hence the condition
\ref{convex criterion} is not necessary.

\subsection{The arcsine distribution on $(-1,1)$}

Recall the arcsine law 
\[
d\mu_{\mathrm{as}}(x)=\frac{\mathbf{1}_{(-1,1)}(x)}{\pi\sqrt{1-x^{2}}}\,dx.
\]
Gross embedding is the unit disk. Thus it is the unique perimeter-minimizing
domain among all bounded rectifiable Jordan domains containing $0$
with real Brownian exit law $\mu_{\mathrm{as}}$. The minimum perimeter
is 
\[
m_{c}(\mu_{\mathrm{as}})=2\pi.
\]
The increasing quantile is 
\[
Q_{\mu_{\mathrm{as}}}(t)=-\cos(\pi t),
\]
so its Gross map is $g(z)=-z$ and its Gross domain is $\D$. Since
$\D$ is convex, theorem \ref{thm:convex-hull-lower-bound} proves
perimeter minimality. However, we give a separate coefficient argument
to establish uniqueness. Let 
\[
f_{U}(z)=\sum_{n\geq1}c_{n}z^{n}
\]
be a conformal map onto an admissible competitor $U$, with $f(0)=0$.
We have 
\[
\frac{1}{2}\sum_{n\geq1}|c_{n}|^{2}=\operatorname{Var}(\mu_{as})=\frac{1}{2},
\]
and therefore 
\[
\sum_{n\geq1}|c_{n}|^{2}=1.
\]
The conformal area formula implies 
\[
\Area(f(\D))=\pi\sum_{n\geq1}n|c_{n}|^{2}\geq\pi.
\]
By the isoperimetric inequality, 
\[
\Per_{c}(f(\D))^{2}\geq4\pi\Area(f(\D))\geq4\pi^{2}.
\]
If $\Per_{c}(f(\D))=2\pi$, equality holds throughout, and hence 
\[
\sum_{n\geq1}(n-1)|c_{n}|^{2}=0.
\]
Thus $c_{n}=0$ for every $n\geq2$ and $|c_{1}|=1$. The image domain
is therefore exactly $\D$.

\subsection{Infinite Skorokhod energy does not imply degeneracy of \ref{optimization problem}.}

In virtue of Proposition \ref{prop:perimeter-bounds}, we saw that
if the Skorokhod energy of $\mu$ is finite then \ref{optimization problem}
is nondegenerate. However, such a condition is not necessary. Let
\[
\Xi(z)=1-\sqrt{1-z},\qquad z\in\mathbb{D}.
\]
We can check that 
\[
\Xi'\in\mathbf{H}^{1}\setminus\mathbf{H}^{2}.
\]
Let $\mu_{\Xi}$ be the law of $\operatorname{Re}\Xi(e^{\Theta i})$
where $\Theta$ is uniform in $(-\pi,\pi)$. On the other hand, $\Xi(\mathbb{D})$
is exactly the Gross embedding of $\mu_{\Xi}$ in virtue of the uniqueness
criterion in \cite{mgremark}. In particular 
\[
\Lambda_{\mu_{\Xi}}=\frac{1}{4}\Vert\Xi'\Vert^{2}_{\mathbf{H}^{2}}=+\infty.
\]
We can also show that $\Xi(\mathbb{D})$ is convex and 
\[
m_{c}(\mu_{\Xi})=\Per_{c}(\Xi(\mathbb{D}))={\textstyle \frac{\sqrt{2}}{2}}\mathbf{Beta}({\textstyle \frac{1}{4}},{\textstyle \frac{1}{2}})\approx3.70815
\]
where $\mathbf{Beta}$ is Euler's beta function.

\begin{figure}[H]
\centering \begin{tikzpicture}[scale=3]
  \fill[blue!15,draw=blue!70!black,thick]
    plot[domain=-45:45,samples=300,variable=\t]
      ({1-sqrt(max(0,2*cos(2*\t)))*cos(\t)},
       {-sqrt(max(0,2*cos(2*\t)))*sin(\t)}) -- cycle;

  \draw[->] (-0.65,0) -- (1.2,0) node[right] {$\operatorname{Re}$};
  \draw[->] (0,-0.65) -- (0,0.65) node[above] {$\operatorname{Im}$};

  \fill (0,0) circle (0.5pt) node[below left,font=\scriptsize] {$\Xi(0)=0$};
  \fill (1,0) circle (0.5pt) node[below right,font=\scriptsize] {$\Xi(1)=1$};
  \pgfmathsetmacro{\xm}{1-sqrt(2)}
  \fill (\xm,0) circle (0.5pt)
    node[below left,font=\scriptsize] {$\Xi(-1)=1-\sqrt2$};

  \node[blue!70!black] at (0.35,0.28) {$\Xi(\mathbb{D})$ };
\end{tikzpicture}

\caption{Gross embedding of $\mu_{\Xi}$, which has the least perimeter. }
\end{figure}
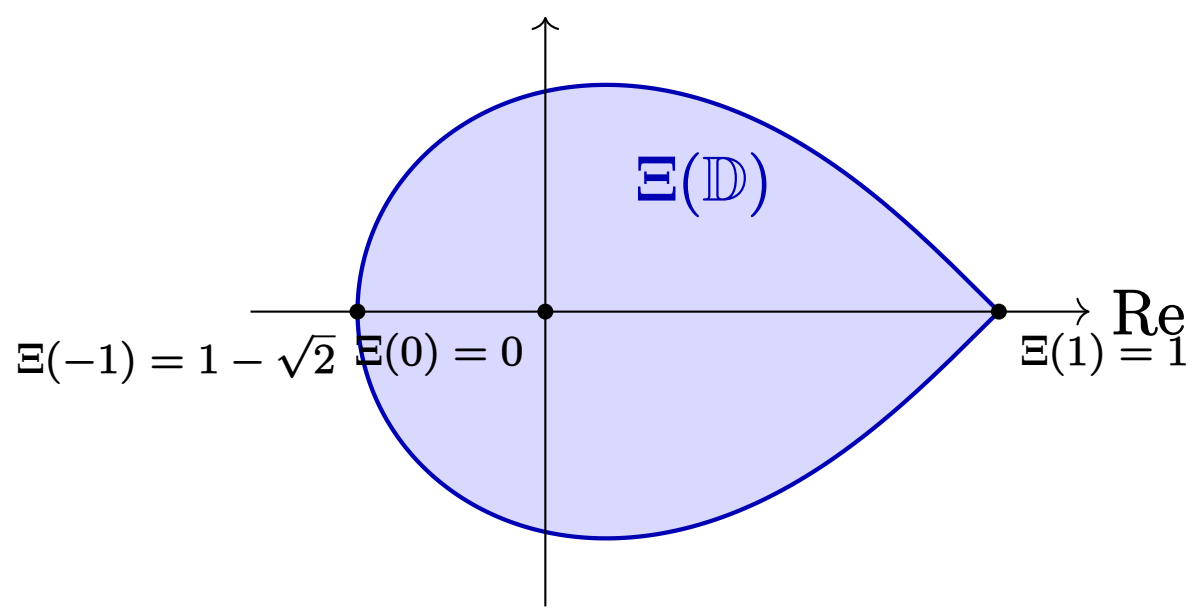

\subsection{Boundedness of the Gross domain does not imply that the perimeter
optimization problem has finite value.}

The classical Cantor-Lebesgue function $C:[0,1]\to[0,1]$ is the continuous,
nondecreasing function that is constant on each complementary interval
of the middle-third Cantor set, satisfies $C(0)=0$ and $C(1)=1$,
and has derivative zero almost everywhere. Define 
\[
Q(t)=t+C(t)-1,\qquad0\le t\le1.
\]
Let $T$ be uniformly distributed on $[0,1]$ and set 
\[
\mu=\operatorname{Law}(Q(T)).
\]
We claim that $\mu$ is centered, absolutely continuous with bounded
density, and has support $[-1,1]$, while 
\[
U^{G}_{\mu}\text{ is bounded},\qquad m_{c}(\mu)=+\infty.
\]

\medskip{}
 \noindent\textit{Properties of the distribution.} Since $C$ is
continuous and nondecreasing, $Q$ is continuous and strictly increasing,
with 
\[
Q(0)=-1,\qquad Q(1)=1.
\]
Moreover, using the symmetry of the Cantor function, 
\[
\int^{1}_{0}Q(t)\,dt=\frac{1}{2}+\frac{1}{2}-1=0.
\]
Thus $\mu$ is centered and $\operatorname{supp}\mu=[-1,1]$.

For $0\le s<t\le1$, we have 
\[
Q(t)-Q(s)=(t-s)+(C(t)-C(s))\ge t-s.
\]
Consequently, $Q^{-1}$ is $1$-Lipschitz. Since the distribution
function of $\mu$ on $[-1,1]$ is precisely $Q^{-1}$, the measure
$\mu$ is absolutely continuous: 
\[
\mu(dx)=\rho(x)\,dx,\qquad0\le\rho(x)\le1\quad\text{a.e.}
\]

\medskip{}
 \noindent\textit{Boundedness of the Gross embedding.} The Cantor
function is Hölder continuous with exponent $\alpha=\log2/\log3$.
Hence $Q$ has the same regularity.

The symmetric decreasing boundary datum associated with the Gross
construction is 
\[
q(\theta)=Q\left(1-\frac{|\theta|}{\pi}\right),\qquad-\pi\le\theta\le\pi.
\]
Thus $q\in C^{\alpha}(\mathbb{T})$. Since the periodic Hilbert transform
preserves Hölder spaces $C^{\alpha}$, $0<\alpha<1$, its harmonic
conjugate is bounded. It follows that the Gross conformal map is bounded
as well, and therefore $G_{\mu}$ is bounded.

\medskip{}
 \noindent\textit{Degeneracy of the perimeter problem.} Let $K\subset[0,1]$
denote the middle-third Cantor set and put $E=Q(K)$. Since $\mathcal{L}^{1}(K)=0$
and $Q$ is injective, 
\[
\mu(E)=\mathcal{L}^{1}(K)=0.
\]

On the other hand, $\mathcal{L}^{1}(E)=1$. Indeed, the complementary
intervals of $K$ have total length $1$, and $C$ is constant on
each of these intervals. Hence $Q$ preserves their lengths. Since
$Q([0,1])=[-1,1]$, whose length is $2$, we obtain 
\[
\mathcal{L}^{1}(Q(K))=2-1=1.
\]

Thus $E\subset\operatorname{supp}\mu$ satisfies 
\[
\mu(E)=0,\qquad\mathcal{L}^{1}(E)>0.
\]
Therefore, 
\[
\mathcal{L}^{1}|_{\operatorname{supp}\mu}\not\ll\mu.
\]
By Theorem~\ref{thm:necessary-density}, no admissible embedding
can have finite conformal perimeter. Consequently, 
\[
m_{c}(\mu)=+\infty.
\]

\begin{figure}[H]
\begin{centering}
\includegraphics[width=8cm,totalheight=8cm,keepaspectratio]{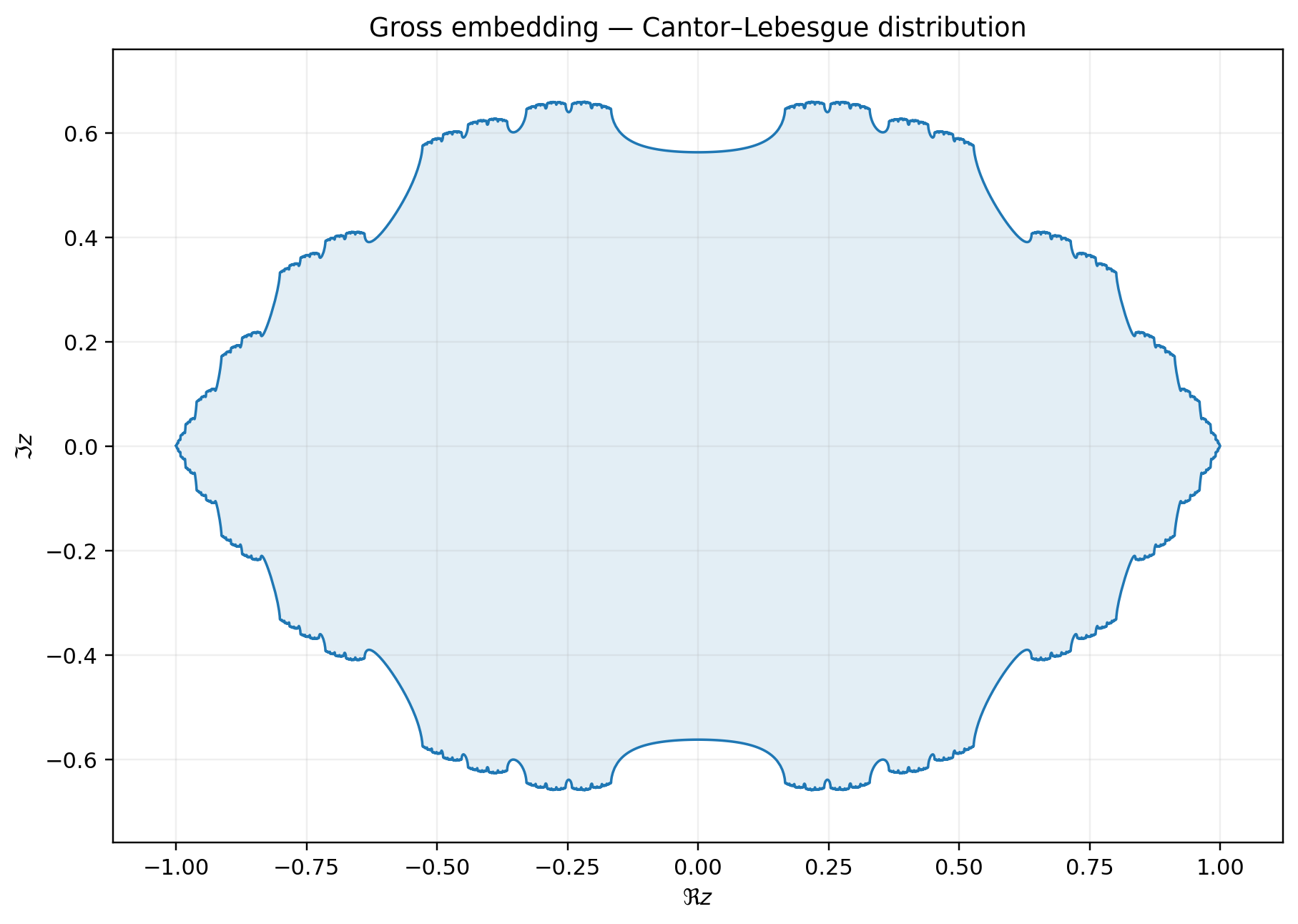} 
\par\end{centering}
\caption{Gross embedding.}
\end{figure}

\section{Concluding remarks and open problems }\label{conclusion}

The planar Skorokhod embedding problem provides a natural setting
in which to investigate how a prescribed probability law constrains
the geometry of a domain. From this perspective, conformal perimeter
is not merely a boundary functional: it measures the geometric cost
of realizing the prescribed real-part exit distribution. The results
of this paper show that this cost is governed by an interplay between
conformal mapping, rearrangement inequalities, and the measure-theoretic
structure of the exit law.

Our findings clarify the role of the Gross embedding. Although it
is distinguished by its area-minimizing property, our explicit counterexamples
show that it need not minimize conformal perimeter. Nevertheless,
it remains central to the optimization problem: its convex hull provides
a universal geometric lower bound, its unboundedness rules out every
finite-cost competitor, and, when it is bounded and convex, it is
the unique perimeter minimizer. Thus, the Gross construction retains
significant extremal information even when it does not itself solve
the perimeter problem.\\

The existence theorem supplies a core result: once a single finite-conformal-perimeter
embedding exists, a minimizer exists as well. No Jordan-boundary assumption
is needed, and every finite-cost minimizer is bounded and admits an
absolutely continuous conformal boundary parametrization. The principal
issue is therefore not attainment in the finite-cost regime, but identifying
which probability laws belong to that regime. Our necessary condition
\[
\mathcal{L}^{1}|_{\operatorname{supp}\mu}\ll\mu
\]
and our logarithmic sufficient condition provide two different descriptions
of this question. The examples demonstrate that neither connected
support nor absolute continuity, even together with boundedness of
the Gross domain, suffices. Conversely, logarithmic integrability
allows finite-perimeter embeddings for densities with substantial
degeneracies, well beyond the class covered by finite Skorokhod energy.
These results leave several natural questions open.

\paragraph{Preservation of finite conformal perimeter.}

The central unresolved question is whether 
\[
m_{\mathrm{c}}(\mu)<\infty\quad\Longrightarrow\quad\operatorname{Per}_{\mathrm{c}}(U^{G}_{\mu})<\infty.
\]
The reverse implication is immediate. A positive answer would therefore
characterize finite-cost solvability through the canonical Gross embedding,
without asserting that Gross minimizes perimeter. A negative answer
would exhibit a different phenomenon: Gross symmetrization could transform
an embedding of finite conformal perimeter into one of infinite conformal
perimeter while preserving the exit law.

Analytically, the difficulty lies at the Hardy-space endpoint. If
$u=\operatorname{Re}f^{*}_{U}$ and $u^{\#}$ is its symmetric decreasing
rearrangement, finite conformal perimeter gives 
\[
u'\in L^{1}(\mathbb{T}),\qquad\mathcal{H}u'\in L^{1}(\mathbb{T}).
\]
Rearrangement controls $(u^{\#})'$ in $L^{1}$, but the additional
conclusion 
\[
\mathcal{H}\bigl((u^{\#})'\bigr)\in L^{1}(\mathbb{T})
\]
does not follow from the usual $L^{1}$ estimates. Understanding whether
the conformal and equimeasurability constraints supply the missing
control would substantially advance the characterization of finite-cost
exit laws.

\paragraph{Optimality beyond the convex case.}

Bounded convex Gross domains are uniquely perimeter minimizing. Is
convexity also necessary for finite-cost perimeter minimality of the
Gross embedding? More generally, can one characterize those probability
laws for which 
\[
m_{\mathrm{c}}(\mu)=\operatorname{Per}_{\mathrm{c}}(G_{\mu})?
\]
Our counterexamples show that this equality does not hold universally,
but they do not determine the full class of laws for which it does.
A characterization would distinguish the distributions for which the
canonical rearrangement is geometrically optimal from those requiring
a different boundary arrangement.

\paragraph{Geometry and uniqueness of minimizers.}

Outside the convex-Gross setting, the geometry of minimizers remains
largely undetermined. Must at least one minimizer be vertically convex?
Must every minimizer have a Jordan boundary, or can a minimizing conformal
boundary parametrization necessarily traverse some boundary portions
more than once? What additional conditions on the exit law imply stronger
boundary regularity?

Uniqueness is closely linked to these questions. Reflection about
the real axis preserves both the exit law and conformal perimeter,
so a unique minimizer must possess this symmetry. If it is also vertically
convex, it is Steiner symmetric and hence coincides with the Gross
embedding. Consequently, for a law whose Gross embedding is not minimizing,
the existence of a vertically convex minimizer would imply nonuniqueness.
This connects the geometry of optimal domains with the structure of
the entire minimizing family.

\paragraph{Quantitative estimates, stability, and approximation.}

The geometric lower bounds, the Skorokhod-energy upper bound, and
the logarithmic sufficient condition provide initial quantitative
control of $m_{\mathrm{c}}(\mu)$. It remains to determine their sharpness
and to understand how the optimal value depends on the prescribed
law. Under what assumptions does convergence of probability measures
imply convergence of their optimal perimeters? When minimizers are
unique, can one obtain corresponding convergence of their normalized
conformal maps? Does the $p$-Brownian convergence concept, developed
in \cite{Brownian conv}, hold? 

These stability questions are also relevant to the construction of
rigorous numerical methods. A useful approximation scheme must control
the exit-law constraint and conformal perimeter while preserving univalence.
Establishing convergence to a genuine minimizer, rather than merely
producing low-perimeter approximate embeddings, would turn the existence
theory into a constructive approach to optimal design.

Taken together, these questions suggest that perimeter optimization
is not a direct extension of area minimization. It introduces a distinct
sensitivity to boundary variation, for which the canonical rearrangement
may preserve some geometric features while worsening others. Understanding
this distinction offers a route toward a broader theory of optimal
planar Skorokhod embeddings, connecting the distribution of probability
mass with the regularity, geometry, and multiplicity of its optimal
realizations.

\section{Appendix: Proofs of the technical tools }\label{Appendix:-Proofs-of}
\begin{proof}
of lemma \ref{lem:hardy} Write $h(z)=H(z^{2})$. A change of angle
in the defining circle means shows $H\in\mathbf{H}^{1}$ and $\|H\|_{\mathbf{H}^{1}}=\|h\|_{\mathbf{H}^{1}}$.
Apply \eqref{eq:classical-hardy} to $H$ at $r^{2}$. 
\end{proof}

For the zero-free functions used below, there is a direct square-root
proof. If $K^{2}=h$, then $K(-z)/K(z)$ is constant with value $1$
at $0$, so $K$ is even and $K(z)=\sum_{n\ge0}d_{n}z^{2n}$. Consequently
\[
|h(r)|=|K(r)|^{2}\le\left(\sum_{n\ge0}|d_{n}|^{2}\right)\left(\sum_{n\ge0}r^{4n}\right)=\frac{\|h\|_{\mathbf{H}^{1}}}{1-r^{4}}.
\]
The zero-free assumption is needed for this particular square-root
argument, not for the lemma itself. The factor $1-r^{4}$ instead
of $1-r^{2}$ is an elementary symmetry improvement of a classical
estimate, not a new general Hardy inequality. 
\begin{proof}
of lemma \ref{lem:boundary} The function $q$ is the decreasing angular
quantile of the law of $q(\Theta)$ for $\Theta$ uniform on $(0,\pi)$.
This law is centered and nondegenerate, so Gross's theorem applies.
Three integrations by parts in the periodic Fourier formula give $c_{n}=O(n^{-3})$.
Hence the series for $g$ and $g'$ converge uniformly on $\overline{\D}$,
and $\Re g(e^{it})=q(t)$. The coefficients are real, and 
\[
g'(0)=c_{1}=-\frac{2}{\pi}\int^{\pi}_{0}q'(t)\sin t\,dt>0.
\]
Since $g'$ never vanishes in $\D$, it is positive on $(-1,1)$.
The endpoint values are $g(1)=q(0)$ and $g(-1)=q(\pi)$; thus 
\begin{equation}
g((-1,1))=(q(\pi),q(0)).\label{eq:real-diameter}
\end{equation}
For $z$ in the upper half-disc, $g(z)$ cannot be real, since otherwise
$g(z)=g(\overline{z})$ would contradict univalence. Its imaginary
part has a constant sign there, and $g'(0)>0$ makes that sign positive.
If $\zeta\in\T$, continuity implies $g(\zeta)\in\overline{g(\D)}$.
It cannot lie in $g(\D)$: continuity of the conformal inverse at
that point, applied to $g(r\zeta)$, would contradict $r\zeta\to\zeta\notin\D$.
Consequently $g(\zeta)\in\partial g(\D)$. For $0<t<\pi$, the imaginary
part of $g(e^{it})$ cannot be zero, because its real part $q(t)$
lies in the interior interval \eqref{eq:real-diameter}.

The upper boundary arc therefore lies strictly above the real axis
and has a strictly decreasing real coordinate. It is injective. The
lower arc is its reflection, and the two arcs meet only at their distinct
endpoints. This proves injectivity on $\overline{\D}$. Their union
is the whole boundary: any boundary point is the limit of interior
images, and compactness of $\overline{\D}$ supplies a preimage on
$\T$. Finally, the $C^{1}$ boundary parametrization has length 
\[
\int^{2\pi}_{0}\left|\frac{d}{dt}g(e^{it})\right|\,dt=\int^{2\pi}_{0}|g'(e^{it})|\,dt=2\pi\|g'\|_{\mathbf{H}^{1}}.
\]
\end{proof}

We will also use the following derivative formula. Integration by
parts in \eqref{eq:gross-series} gives 
\[
nc_{n}=-\frac{2}{\pi}\int^{\pi}_{0}q'(t)\sin(nt)\,dt.
\]
For $0<r<1$, summing the uniformly convergent geometric series 
\[
\sum_{n\ge1}r^{n-1}\sin(nt)=\frac{\sin t}{1-2r\cos t+r^{2}}
\]
therefore yields 
\begin{equation}
g'(r)=-\frac{2}{\pi}\int^{\pi}_{0}\frac{q'(t)\sin t}{1-2r\cos t+r^{2}}\,dt.\label{eq:derivative-kernel}
\end{equation}

\begin{proof}
of corollary \ref{cor:isoperimetric}. The area inequality \eqref{eq:area-monotone}
gives 
\[
\frac{\Area(W_{a,b})}{\Area(G_{a,b})}\geq1.
\]
Combining this with \eqref{eq:family-upper} and \eqref{eq:family-lower},
we obtain 
\begin{align*}
\frac{I(G_{a,b})}{I(W_{a,b})} & =\left(\frac{\Per(G_{a,b})}{\Per(W_{a,b})}\right)^{2}\frac{\Area(W_{a,b})}{\Area(G_{a,b})}\\
 & \geq\left(\frac{\Per(G_{a,b})}{\Per(W_{a,b})}\right)^{2}\\
 & >\left(\frac{\frac{8}{15}(1+a^{2})(5+3b)}{2\pi(1-b)+8b}\right)^{2}=R(a,b)^{2}.
\end{align*}
This proves \eqref{eq:iso-family}. If $R(a,b)\geq1$, the strict
inequality implies $I(G_{a,b})>I(W_{a,b})$. 
\end{proof}

\begin{proof}
of Theorem \ref{thm:family}.\emph{ The original domain.} Since $ia\in\D$,
the translation in \eqref{eq:family} puts $0$ inside $W_{a,b}$.
If $b<1$, then 
\[
\Re F_{b}'(z)=\Re(1+bz^{2})\ge1-b>0\qquad(z\in\overline{\D}).
\]
For distinct $z,w\in\overline{\D}$, integration along their joining
segment shows that 
\[
\Re\frac{F_{b}(z)-F_{b}(w)}{z-w}=\int^{1}_{0}\Re F_{b}'(w+s(z-w))\,ds>0.
\]
Thus $F_{b}$ is injective on $\overline{\D}$ and its derivative
is nonvanishing there. When $b=1$, equality $F_{1}(z)=F_{1}(w)$
with $z\ne w$ would imply $z^{2}+zw+w^{2}=-3$. Equality in the triangle
inequality then forces $z^{2}=zw=w^{2}=-1$, and hence $z=w$, a contradiction.
Also $F_{1}'$ does not vanish in $\D$. This proves the Jordan and
conformal assertions for the original domains. Translation does not
affect length, and 
\[
|1+be^{2it}|=|(1-b)+b(1+e^{2it})|\le(1-b)+2b|\cos t|.
\]
Integration proves \eqref{eq:family-upper}, with equality $8$ when
$b=1$.

\medskip{}
 \emph{The exit law and its quantile.} The map $F_{b}-F_{b}(ia)$
sends $ia$ to $0$. Conformal invariance of planar Brownian motion
therefore identifies the exit law from $W_{a,b}$ with the pushforward
of harmonic measure at $ia$ in the disc; see \cite{MP2010,thesis,greg notes}.
The corresponding angular density is the Poisson kernel 
\[
p_{a}(t)=\frac{1-a^{2}}{2\pi(1+a^{2}-2a\sin t)}.
\]
The real boundary coordinate is 
\begin{equation}
r_{b}(t)=\Re F_{b}(e^{it})=\cos t+\frac{b}{3}\cos(3t)=(1-b)\cos t+\frac{4b}{3}\cos^{3}t.\label{eq:real-boundary}
\end{equation}
It is even and strictly decreasing on $[0,\pi]$. Indeed, 
\[
r_{b}'(t)=-\sin t\bigl(1-b+4b\cos^{2}t\bigr),
\]
which is negative on $(0,\pi)$ for $b<1$; for $b=1$ the function
$\frac{4}{3}\cos^{3}t$ is still strictly decreasing. Folding the
angles $t$ and $-t$ gives, for $0<t<\pi$, 
\begin{equation}
p_{a}(t)+p_{a}(-t)=\frac{k}{\pi(\cos^{2}t+k^{2}\sin^{2}t)}.\label{eq:folded-density}
\end{equation}
Define the increasing continuous angle $\theta=\theta(t)$ on $[0,\pi]$
by 
\begin{equation}
\tan\theta=k\tan t,\qquad\theta(0)=0,\quad\theta(\pi/2)=\pi/2.\label{eq:angle-change}
\end{equation}
The continuous branch is used across $\pi/2$. Since 
\[
\theta'(t)=\frac{k}{\cos^{2}t+k^{2}\sin^{2}t},
\]
\eqref{eq:folded-density} becomes $d\theta/\pi$. Therefore the decreasing
angular quantile is 
\begin{equation}
q_{a,b}(\theta)=(1-b)C_{k}(\theta)+\frac{4b}{3}C_{k}(\theta)^{3},\qquad C_{k}(\theta)=\frac{\cos\theta}{\sqrt{\cos^{2}\theta+k^{-2}\sin^{2}\theta}}.\label{eq:family-quantile}
\end{equation}
Its even periodic extension is real-analytic (in particular $C^{3}$),
strictly decreasing on $[0,\pi]$, and satisfies 
\begin{equation}
q_{a,b}(\pi-\theta)=-q_{a,b}(\theta).\label{eq:quantile-symmetry}
\end{equation}
The underlying law is consequently centered and symmetric, with support
$[-1-b/3,1+b/3]$.

Let $g=g_{a,b}$ be the analytic completion of $q=q_{a,b}$. Lemma~\ref{lem:boundary}
shows that $G_{a,b}=g(\D)$ is a bounded rectifiable Jordan domain.
By construction, its real exit law agrees with that of $W_{a,b}$.
Equation~\eqref{eq:quantile-symmetry} makes every even-indexed coefficient
of $g$ vanish. Thus $g$ is odd and $g'$ is even and zero-free.

\medskip{}
 \emph{A lower bound for one derivative value.} Pairing $t$ and $\pi-t$
in \eqref{eq:derivative-kernel}, and using $q'(\pi-t)=q'(t)$, gives
\begin{equation}
g'(a)=-\frac{4(1+a^{2})}{\pi}\int^{\pi/2}_{0}\frac{q'(t)\sin t}{(1-2a\cos t+a^{2})(1+2a\cos t+a^{2})}\,dt.\label{eq:paired-kernel}
\end{equation}
Use the substitution $t=\arctan(kx)$. Then 
\begin{align*}
q(\arctan(kx)) & =\frac{1-b}{\sqrt{1+x^{2}}}+\frac{4b}{3(1+x^{2})^{3/2}},\\
q'(t)\,dt & =-x\left[\frac{1-b}{(1+x^{2})^{3/2}}+\frac{4b}{(1+x^{2})^{5/2}}\right]dx,\\
\sin t & =\frac{kx}{\sqrt{1+k^{2}x^{2}}},\\
(1-2a\cos t+a^{2})(1+2a\cos t+a^{2}) & =(1-a^{2})^{2}\cos^{2}t+(1+a^{2})^{2}\sin^{2}t\\
 & =(1-a^{2})^{2}\frac{1+x^{2}}{1+k^{2}x^{2}}.
\end{align*}
Substitution into \eqref{eq:paired-kernel}, together with $(1+a^{2})k=1-a^{2}$,
yields the exact identity 
\begin{equation}
\begin{split}g'(a)=\frac{4}{\pi(1-a^{2})}\int^{\infty}_{0}x^{2}\sqrt{1+k^{2}x^{2}}\left[\frac{1-b}{(1+x^{2})^{5/2}}+\frac{4b}{(1+x^{2})^{7/2}}\right]dx.\end{split}
\label{eq:family-derivative}
\end{equation}
All terms are nonnegative, and $\sqrt{1+k^{2}x^{2}}>1$ for $x>0$.
The substitution $u=x/\sqrt{1+x^{2}}$ gives 
\begin{align*}
\int^{\infty}_{0}\frac{x^{2}}{(1+x^{2})^{5/2}}\,dx & =\int^{1}_{0}u^{2}\,du=\frac{1}{3},\\
\int^{\infty}_{0}\frac{x^{2}}{(1+x^{2})^{7/2}}\,dx & =\int^{1}_{0}u^{2}(1-u^{2})\,du=\frac{2}{15}.
\end{align*}
It follows that 
\begin{equation}
g'(a)>\frac{4}{\pi(1-a^{2})}\left(\frac{1-b}{3}+\frac{8b}{15}\right)=\frac{4(5+3b)}{15\pi(1-a^{2})}.\label{eq:family-derivative-lower}
\end{equation}
Apply \eqref{eq:perimeter-point} with $r=a$ to conclude 
\[
\Per(G_{a,b})\ge2\pi(1-a^{4})g'(a)>\frac{8}{15}(1+a^{2})(5+3b),
\]
which proves \eqref{eq:family-lower}.

\medskip{}
 \emph{Analytic boundaries when $b<1$.} The assertion for $W_{a,b}$
was established above. The real-analytic periodic function $q$ has
exponentially decreasing Fourier coefficients, so $g$ extends holomorphically
across $\T$. For $0<t<\pi$, strict negativity of $q'(t)$ implies
\[
\Re\bigl(ie^{it}g'(e^{it})\bigr)=q'(t)<0;
\]
hence the boundary derivative does not vanish on the upper arc, or
on the lower arc by reflection. At $1$, let $r\uparrow1$ in \eqref{eq:derivative-kernel}.
Since $q'(t)=O(\sin t)$ at the endpoints, dominated convergence is
applicable and gives 
\[
g'(1)=-\frac{1}{\pi}\int^{\pi}_{0}q'(t)\cot(t/2)\,dt>0.
\]
For example, for $r\ge1/2$ the denominator in \eqref{eq:derivative-kernel}
is at least $1-\cos t$, which supplies a bounded majorant near $0$.
The evenness of $g'$ gives $g'(-1)=g'(1)>0$. Thus $g'$ is nonvanishing
on the entire circle. Together with boundary injectivity, this proves
that $\partial G_{a,b}$ is a regular real-analytic Jordan curve. 
\end{proof}

\end{document}